\documentclass[pdflatex,sn-mathphys-num]{sn-jnl}

\usepackage{graphicx}
\usepackage{amsmath,amssymb,amsfonts}
\usepackage{amsthm}
\usepackage{mathrsfs}
\usepackage{enumitem}
\usepackage[title]{appendix}
\usepackage[dvipsnames]{xcolor}
\usepackage{manyfoot}
\usepackage{booktabs}
\usepackage{listings}

\theoremstyle{thmstyleone}
\newtheorem{theorem}{Theorem}
\newtheorem{lemma}[theorem]{Lemma}
\newtheorem{corollary}[theorem]{Corollary}
\newtheorem{proposition}[theorem]{Proposition} 

\theoremstyle{thmstyletwo}

\newtheorem{remark}{Remark}

\theoremstyle{thmstylethree}

\DeclareMathOperator*{\argmin}{arg\,min}

\DeclareMathOperator{\dist}{dist}
\newcommand{\qbox}[1]{\quad\hbox{#1}\quad}

\newcommand{\R}{\mathbb{R}}

\newcommand{\ackname}{Acknowledgements}
\def\acknowledgement{\par\addvspace{17pt}\small\rmfamily
\trivlist\if!\ackname!\item[]\else
\item[\hskip\labelsep
{\bfseries\ackname}]\fi}

\newenvironment{acknowledgements}{\begin{acknowledgement}}
{\end{acknowledgement}}

\begin{document}

\title{A unified analysis on the speedup of accelerated gradient methods: An inertial dynamics approach}

\author{\fnm{Zepeng} \sur{Wang}}\email{zepeng.wang@rug.nl}

\author{\fnm{Juan} \sur{Peypouquet}}\email{j.g.peypouquet@rug.nl}

\affil{\orgdiv{Bernoulli Institute for Mathematics, Computer Science and Artificial Intelligence}, \orgname{University of Groningen}, \state{Groningen}, \country{The Netherlands}}

\abstract{Nesterov's Accelerated Gradient Method is one of the most popular first order optimization algorithms. When applied to $L$-smooth, $\mu$-strongly convex functions, it converges at a rate of $\mathcal{O}\left(\left(1-\sqrt{\frac{\mu}{L}}\right)^{k}\right)$. The more recent {\it Triple Momentum Method} and the {\it Information Theoretic Exact Method} enjoy an improved rate of
$\mathcal{O}\left(\left(1-2\sqrt{\frac{\mu}{L}}\right)^{k}\right)$. Their analysis relies on {\it integral quadratic constraints} and {\it performance estimation techniques}, respectively.  In this work, we provide a dynamic explanation for this {\it factor 2 speedup} based on the subtle relationships between the coefficients of an inertial system with Hessian-driven damping. The standard explicit discretization of this second order ordinary differential equation produces intuitive variants of Nesterov's method with sped-up convergence rates. The proof strategy allows us to extend the analysis beyond the strongly convex setting to account for convex functions with quadratic growth or the Polyak-{\L}ojasiewicz inequality. Under uniqueness of the minimizer, we establish a factor $\sqrt{2}$ speedup with respect to the state of the art. With no assumption on the set of minimizers, the new convergence rate (asymptotically) matches that of Gradient Descent, a fact that was previously unknown.}



\maketitle

\section{Introduction}
Let $H$ be a real Hilbert space, and consider the convex optimization problem:
\begin{equation*}
\min_{x\in H} f(x),
\end{equation*}
where $f:H\to\R$ is convex and $L$-smooth, which means that $\nabla f$ is Lipschitz-continuous with constant $L>0$. We assume that the solution set $S:\argmin(f)$ is nonempty, and denote $x^*\in S$ and $f^*=\min(f)$. \\

{\bf Nesterov's method for strongly convex functions.} If $f$ is $\mu$-strongly convex, this problem can be solved by Nesterov's {\it Accelerated Gradient Method} (AGM) \cite{Nesterov_2004}:
\begin{equation}\label{E: NAG-SC}
\left\{
\begin{array}{rcl}
y_{k+1} &=& x_k - s\nabla f(x_k),\\[3pt]
x_{k+1} &=& y_{k+1} + \frac{1-\sqrt{q}}{1+\sqrt{q}}(y_{k+1}-y_k),
\end{array}
\right.
\end{equation}
which guarantees an accelerated linear convergence rate \cite{Nesterov_2004,Bansal_2019,Luo_2022}  
$$ f(y_{k+1}) - f^* = \mathcal{O}\left( (1-\sqrt{q})^k \right),\qbox{with} s=1/L \qbox{and} q:=\mu/L.  $$
Siegel \cite{Siegel_2019} obtained a similar convergence rate for a slightly different algorithm whose iterations have four substeps, instead of two. In \cite{Aujol_2022}, Siegel's method was improved to achieve a convergence rate of $\mathcal{O}\left( (1+\sqrt{2q})^{-k}\right)$, which was also obtained in \cite{Luo_2021} for an algorithm composed of three steps and based on a Gauss-Seidel type discretization of an inertial system, and in \cite{Park_2023} for a variant of the so-called {\it Optimized Gradient Method for Strongly Convex functions} (OGM-SC) \cite{Drori_2014,Kim_2016,Kim_2017}, borrowing tools from linear coupling analysis in \cite{Zhu_2017}. The best convergence rate known to date is $\mathcal{O}\left( (1-2\sqrt{q})^{k} \right)$, achieved by the {\it Triple Momentum Method} (TMM) of \cite{Scoy_2018}, whose design and analysis rely on integral quadratic constraints \cite{Lessard_2016}, and by the {\it Information Theoretic Exact Method} (ITEM) of \cite{Taylor_2023_optimal}, based on performance estimation techniques \cite{Drori_2014}. A unified and intuitive explanation for this $\sqrt{q}\to\sqrt{2q}\to 2\sqrt{q}$ speedup has been lacking. \\

{\bf Continuous time models.} The interplay between inertial dynamics and algorithms dates back to Polyak's {\it Heavy Ball Method} \cite{Polyak_1964}, and has received much attention in recent years, especially after \cite{Su_2016} linked Nesterov's method for convex functions \cite{Nesterov_1983} to an inertial system with asymptotic vanishing damping (see also \cite{Jordan_2016,Attouch_2016_Hessian,May_2017,Attouch_2018,Aujol_2019,Siegel_2019,Luo_2022,Shi_2022,Shi_2024}, among others). The convergence analysis of optimization algorithms inspired by the study of continuous-time dynamics is usually intuitive and easy to follow. This is a line of research promoted by Hedy Attouch, to whom we dedicate this work.

Polyak's Heavy Ball Method is obtained by a natural discretization of the {\it Heavy-Ball with Friction} system \cite{Polyak_1964,Attouch_2000}
\begin{equation}\label{E: HBF}
\ddot{x} + 2\sqrt{\mu}\dot{x} + \nabla f(x) = 0.
\end{equation}
It was shown in \cite{Jordan_2021} that this also provides a continuous-time model for
Nesterov's accelerated gradient method \eqref{E: NAG-SC}, although the relationship is less transparent. A more accurate description of \eqref{E: NAG-SC} is given by the {\it High Resolution} dynamical system \cite{Shi_2022}
\begin{equation}\label{E: HBF-H-NAG}
\ddot{x} + 2\sqrt{\mu}\dot{x} + \sqrt{s}\nabla^2 f(x)\dot{x} + \left( 1 + \sqrt{\mu s} \right)\nabla f(x) = 0,
\end{equation}
which incorporates the Hessian-driven damping term
introduced in \cite{Alvarez_2002} with the purpose of taming the heavy oscillations that are typical to the solutions of \eqref{E: HBF} (see also \cite{Attouch_2016_Hessian}). It was shown in \cite{Shi_2022} that the solutions of \eqref{E: HBF-H-NAG} satisfy $f(x(t))-f^*=\mathcal O\left( e^{-\frac{\sqrt{\mu}t}{4}} \right)$. This rate was improved to $\mathcal{O}\left( e^{-\frac{\sqrt{\mu}t}{2}} \right)$ in \cite{Attouch_2022}. However, both are slower than the $\mathcal O\left(e^{-\sqrt{\mu}t}\right)$ obtained in \cite{Siegel_2019} for \eqref{E: HBF}  (see also \cite{Jordan_2021} and \cite{Attouch_2022}). \\

{\bf From continuous dynamics to algorithms.} The high resolution dynamical system \eqref{E: HBF-H-NAG} can be written, more generally, as
\begin{equation}\label{System: HBF-H}\tag{HBF-H}
\ddot{x} + \alpha\dot{x} + \beta\nabla^2 f(x) \dot{x} + \gamma \nabla f(x) = 0,
\end{equation}
with parameters $\alpha,\beta,\gamma > 0$. A finite difference discretization of this ODE with a step size $h>0$ yields
\begin{equation} \label{E:Algorithm_family}
\frac{(x_{k+1}-x_k)-(x_k-x_{k-1})}{h^2} + \alpha\, \frac{x_{k+1}-x_k}{h} + \beta\,\frac{\nabla f(x_k) - \nabla f(x_{k-1})}{h} + \gamma \nabla f(x_k) = 0.
\end{equation}
Although this is not always obvious, Nesterov's AGM \eqref{E: NAG-SC}, as well as the variants in \cite{Siegel_2019,Attouch_2022,Aujol_2022,Park_2023}, can all be cast in this form, with coefficients as shown in the following table:

\begin{center}
\renewcommand{\arraystretch}{1.8}
\begin{tabular}{|c|c|c|c|}
\hline
{\bf Method} & $\alpha$ & $\beta$ & $\gamma$  
 \\
\hline
 {\cite{Nesterov_2004}}  
 & $\alpha h = \frac{2\sqrt{q}}{1-\sqrt{q}}$ & $\beta=h$ & $\gamma = \frac{1+\sqrt{q}}{1-\sqrt{q}}$\\
\hline
 {\cite{Siegel_2019}}  
 & $\alpha h = 2\sqrt{q} + q$ & $\beta=h$ & $\gamma = 1 + \sqrt{q}$ \\
\hline
 {\cite{Aujol_2022}}  
 & $\alpha h = \frac{\frac{3}{2}\sqrt{2q}}{1+2\sqrt{2q}}$ & $\beta=h$ &  $\gamma = \frac{1 + \frac{7}{2}\sqrt{2q}}{1+3\sqrt{2q}}$ \\
\hline
 {\cite{Park_2023}}  
 & $\alpha h = \frac{\sqrt{8q + q^2} + 3q }{1-q}$ & $\beta=h$ &  $\gamma = \frac{ 2 + \sqrt{8q + q^2} + q }{1-q}$ \\
\hline 
 {\cite{Attouch_2022}}    
 & $\alpha h = \frac{2\sqrt{q}}{1-\sqrt{q}}$ & $\beta=\frac{h}{1-\sqrt{q}}$ & $\gamma = \frac{1}{1-\sqrt{q}}$ \\
 \hline
\end{tabular}
\end{center}
With the exception of \cite{Attouch_2022}, all these algorithms use the normalization $\beta=h=\frac{1}{\sqrt{L}}$, whereas the choice of $\alpha$ and $\gamma$ is more diverse. The algorithms in \cite{Scoy_2018,Taylor_2023_optimal}, which have the fastest proved convergence rate $\mathcal{O}\left( (1-2\sqrt{q})^{k} \right)$, seem to be of a different nature. \\

The three coefficients in \eqref{System: HBF-H} can be reduced to two by either a function or a time rescaling (see also \cite{Attouch_2019_proximal,Attouch_2022_rescaling}). Although most works on this kind of dynamical systems normalize $\gamma=1$, we prefer to keep the three-parameter representation for several reasons: First, the hidden coefficient reappears via the step size $h$ when passing to the algorithmic setting. Second, this structure allows for a more direct comparison with other methods in the literature, where the normalization $\beta=h$ is the standard practice. Finally, the freedom to use $\gamma$ as a design variable allows us to identify algorithms within the family described by \eqref{E:Algorithm_family}, which achieve the fastest convergence rate apparently reserved for \cite{Scoy_2018,Taylor_2023_optimal}. It turns out that, as the impact of $\gamma$ on the speedup described above becomes transparent, one can extend the analysis beyond strongly convex functions to observe this phenomenon arise under different geometries. \\

{\bf Our contribution.} The starting point of the present research was to better understand, from the perspective of inertial dynamics, the parameter choices that lead to faster convergence rates. 
\begin{enumerate}
    \item We unify the Lyapunov analysis for \eqref{System: HBF-H}, and explain how the components of the energy function can be tuned in order to explain the speedup. For example, we establish a convergence rate of $f(x)-f^*=\mathcal O\left(e^{-\sqrt{\mu}t}\right)$ for \eqref{E: HBF-H-NAG}, improving on \cite{Shi_2022} by a factor of $4$, and removing the gap with respect to \eqref{E: HBF}. If one fixes $\gamma\sim 1$ (as $\mu\to 0$) to account for rescaling, we improve on existing rates by an additional factor of $\sqrt{2}$. If one increases the value of $\gamma$ from $\gamma\sim 1$ to $\gamma\sim 2$ (as $\mu\to 0$), one can improve the convergence rate by yet an additional factor of $\sqrt{2}$. For $\gamma\sim 2$, we derive a fast convergence rate $f(x)-f^*=\mathcal O\left(e^{-2\sqrt{\mu}t}\right)$. The range $\gamma\in[1,2]$ is important, as we explain later. In the absence of strong convexity, one must accept more conservative convergence rate guarantees \cite{Lin_2023}. Nevertheless, convergence rate analysis of first-order methods under relaxed strong convexity conditions \cite{Aujol_2022,Aujol_2023,Bolte_2017,Zhang_2017,Karimi_2016,Lewis_2018,Necoara_2019} is an interesting topic with relevant consequences in many practical applications. With this in mind, we obtain analogous results for functions with quadratic growth (we refer the reader to Appendix \ref{A:geometry} for further details, including its connection with strong convexity and the Polyak-{\L}ojasiewicz inequality). As before, we recover the state-of-the-art rates and identify the mechanism that makes speedup possible. The estimations for general $\gamma$ bring insight that can be transposed to the algorithmic setting.

    \item Using the insight gained from the continuous-time models, we refine the Lyapunov analysis to derive convergence rates for the algorithms in the family defined by \eqref{E:Algorithm_family}. We identify regions in the parameter space that guarantee the stability of the discretization and the subsequent preservation of the linear convergence. In the strongly convex case, the analysis is valid for $\gamma\in[1,2]$, whence the importance of this interval. By appropriately tuning the parameters, one obtains algorithms with the fastest convergence rates known to date, which asymptotically match \cite{Scoy_2018,Taylor_2023_optimal}. For functions with quadratic growth, we strictly improve upon existing results. For example, for convex functions with a unique minimizer, we obtain a convergence rate of $f(x)-f^*=\mathcal O\left( (1+2(\sqrt{2}-1)\sqrt{q})^{-k} \right)$, which gives a factor $\sqrt{2}$ improvement with respect to the best convergence rate so far, namely $\mathcal{O}\left( (1+(2-\sqrt{2})\sqrt{q})^{-k} \right)$ \cite{Aujol_2023}.
    
\end{enumerate}

{\bf Organization of the paper.} Section \ref{Sec: inertial_dynamics} is concerned with the study of the inertial dynamics \eqref{System: HBF-H}. We define a master Lyapunov function and explain how to particularize it to obtain a factor $\sqrt{2}$ speedup in Subsection \ref{SS:factorsqrt2} and a factor $2$ in \ref{SS:factor2}. Subsections \ref{SS:quadunique} and \ref{SS:quadnonunique} deal with functions with quadratic growth and a unique or multiple minimizers, respectively. Section \ref{Sec: inertial_algorithms} has a similar structure and is devoted to the analysis of the Nesterov-type algorithms generated via \eqref{E:Algorithm_family}.

\section{Inertial dynamics}\label{Sec: inertial_dynamics}
In this section, we prove linear convergence rates of the function values along the trajectories of system \eqref{System: HBF-H}. Throughout this section, we assume that $f:H\to\R$ is twice differentiable, and that the trajectory $x:[0,\infty)\to H$ satisfies \eqref{System: HBF-H} with initial conditions $x(0)=x_0$ and $\dot{x}(0)=v_0$. 

As it is unlikely to lead to confusion, we omit the explicit dependence on $t$ of the time-dependent variables hereafter to alleviate the notation. For example, we write $x$ for $x(t)$ and $f(x)$ for $f\left( x(t) \right)$.

Our analysis relies on the energy function $W:[0,\infty)\to\R$, defined by 
\begin{equation}\label{E: W}
W = \frac{1}{2} \| \dot{x} + \beta\nabla f(x) + \xi( x - x^* ) \|^2 - \frac{\eta}{2} \| x - x^* \|^2 + \theta ( f(x) - f^* ),
\end{equation}
where $0\le \xi\le \alpha$, $\eta = \omega\xi(\alpha-\xi)$, $0\le\omega\le 1$ and $\theta = \gamma - (\alpha-\xi)\beta> 0$. The values of $\xi$ and $\omega$ will be fixed later. 

\begin{remark}\label{Rem: func_bound}
If $f$ satisfies $f(x) - f^*\ge \frac{\mu}{2}\|x-x^*\|^2$ for every $x\in H$,\footnote{For further insight into this inequality, the reader is referred to Appendix \ref{A:geometry}.} then
$$- \frac{\eta}{2} \| x - x^* \|^2 + \theta ( f(x) - f^* ) \ge \left( \theta - \frac{\eta}{\mu} \right)( f(x) - f^* ). $$
As a consequence,
$$ f\left(x(t)\right) - f^* \le \frac{\mu}{\mu\theta- \eta}W(t), $$
provided $\eta<\theta\mu$.
\end{remark}

The following result will facilitate the subsequent convergence analysis. The proof is technical, and given in Appendix \ref{A:Lemma_continuous}.

\begin{lemma}\label{Lem: W_bound}
Let $W$ be defined by \eqref{E: W}. Then,
\begin{align*}
\dot{W} + (1+\omega)(\alpha-\xi)W &\le - \beta\left( \theta - \frac{\omega}{2} \right)\|\nabla f(x)\|^2 
  + (1+\omega)(\alpha-\xi)\theta (f(x)-f^*)\\
&\quad 
  - \xi\theta \langle\nabla f(x),x-x^* \rangle \\
&\quad + \frac{1}{2}(\alpha-\xi)\xi[ (1+\omega)(\xi-\omega(\alpha-\xi))  + \omega\xi(\alpha-\xi)\beta ]\|x-x^*\|^2.
\end{align*}
\end{lemma}

\subsection{ Strong convexity: factor \texorpdfstring{$\sqrt{2}$}{} acceleration }\label{SS:factorsqrt2}
In this subsection, we will show how to accelerate the convergence of system \eqref{E: HBF-H-NAG} by a factor of $\sqrt{2}$, namely: from $\mathcal{O}\left( e^{-\sqrt{\mu}t} \right)$ to $\mathcal{O}\left( e^{-\sqrt{2\mu}t} \right)$. The key is to use the function $W$, defined by \eqref{E: W}, with $\xi=\frac{\alpha}{2}$ and $\omega=0$. \\

The main result of this subsection is the following:

\begin{theorem}
Let $f:H\to\mathbb{R}$ be $\mu$-strongly convex. Pick $\theta>0$, and consider the system \eqref{System: HBF-H}, where
$$\alpha\in \left(0,\, 2\sqrt{\mu\theta}\right],\quad \beta > 0,\quad \gamma = \theta + \frac{1}{2}\alpha\beta.$$
If $v_0 = -\beta\nabla f(x_0)$, we have, for every $t\ge 0$,
$$f\left( x(t) \right) - f^* \le 2 \left( f(x_0) - f^* \right) e^{-\frac{1}{2}\alpha t}.$$
Moreover, 
$$ \int_0^{\infty} e^{\frac{\alpha}{2}t}\| \nabla f\left(x(t)\right) \|^2 dt < \infty. $$
\end{theorem}

\begin{proof}
Let $W_1$ be the realization of the function $W$, defined by \eqref{E: W}, with $\xi=\frac{\alpha}{2}$ and $\omega=0$, namely:
$$W_1(t) = \frac{1}{2} \left\| \dot{x} + \beta\nabla f(x) + \frac{\alpha}{2}( x - x^* ) \right\|^2 + \theta ( f(x) - f^* ).$$
From Lemma \ref{Lem: W_bound}, we obtain
\begin{align*}
\dot{W_1} + \frac{\alpha}{2}W_1 
&\le -\beta\theta\| \nabla f(x) \|^2 
  + \frac{\alpha\theta}{2} (f(x)-f^*)
  - \frac{\alpha\theta}{2} \langle\nabla f(x),x-x^* \rangle 
  + \frac{\alpha^3}{16}\|x-x^*\|^2.
\end{align*} 
Using strong convexity:
$$ \langle\nabla f(x),x-x^* \rangle \ge f(x)-f^* + \frac{\mu}{2}\|x-x^*\|^2,$$
and with $\alpha\le2\sqrt{\mu\theta}$ in mind, we arrive at
$$ \dot{W_1} + \frac{\alpha}{2}W_1 + \beta\theta\| \nabla f(x) \|^2 
\le \frac{\alpha (\alpha^2 - 4\mu\theta)}{16}\|x-x^*\|^2
\le 0. $$
Since $\beta,\theta>0$, we have 
$$ W_1(t)\le W_1(0)e^{-\frac{1}{2}\alpha t} \qbox{and} \int_0^{\infty} e^{\frac{\alpha}{2}t}\| \nabla f\left(x(t)\right) \|^2 dt < \infty. $$
Using $v_0 = -\beta\nabla f(x_0)$ , we obtain
\begin{align*}
W_1(0) 
&= \frac{\alpha^2}{8}\| x_0 - x^* \|^2 + \theta \left( f(x_0) - f^* \right)\\
&\le \frac{\mu\theta}{2}\| x_0 - x^* \|^2 + \theta \left( f(x_0) - f^* \right) \\
&\le 2\theta \left( f(x_0) - f^* \right),
\end{align*}
in view of $\alpha\le2\sqrt{\mu\theta}$ and the strong convexity. It follows that
$$ f\left(x(t) \right) - f^* \le \frac{W_1(t)}{\theta} \le \frac{W_1(0)}{\theta}e^{-\frac{1}{2}\alpha t} \le 2 \left( f(x_0) - f^* \right) e^{-\frac{1}{2}\alpha t}, $$
as claimed.
\end{proof}

Setting $\theta=1$ and $\alpha = 2\sqrt{\mu}$, we  immediately obtain an improved convergence rate for the high resolution dynamical system \eqref{E: HBF-H-NAG} by a factor of $4$ with respect to \cite{Shi_2022}). This matches the rate obtained in \cite[Theorem 1]{Siegel_2019}, \cite[Proposition 3]{Jordan_2021} and \cite[Theorem 7]{Attouch_2022} for the heavy ball with friction system \eqref{E: HBF}, closing the existing gap in the theory:

\begin{corollary}
Let $f:H\to\mathbb{R}$ be $\mu$-strongly convex, and let $x$ be a solution of \eqref{E: HBF-H-NAG} with $v_0 = -\sqrt{s}\nabla f(x_0)$. For every $t\ge 0$, we have
$$f\left( x(t) \right) - f^* \le 2 \left( f(x_0) - f^* \right) e^{-\sqrt{\mu}t}.$$
Moreover, 
$$ \int_0^{\infty} e^{\sqrt{\mu}t}\| \nabla f\left(x(t)\right) \|^2 dt < \infty. $$    
\end{corollary}

\begin{remark} 
For $\theta = 2$ and $\alpha=2\sqrt{2\mu}$, one obtains $f(x) - f^* = \mathcal{O}\left( e^{-\sqrt{2\mu} t} \right)$ for the system 
$$\ddot{x} + 2\sqrt{2\mu}\dot{x} + \beta\nabla^2 f(x) \dot{x} + \left( 2 + \sqrt{2\mu}\beta \right) \nabla f(x) = 0.$$    
\end{remark}

\subsection{ Strong convexity: factor \texorpdfstring{$2$}{} acceleration } \label{SS:factor2}
In this subsection, we will show how to attain the convergence rate $\mathcal{O}\left( e^{-2\sqrt{\mu} t} \right)$ for system \eqref{System: HBF-H}, which holds the key to establishing the state-of-the-art rate of TMM \cite{Scoy_2018} and ITEM \cite{Taylor_2023_optimal}. To this end, we use the function $W$, defined by \eqref{E: W}, but now with $\xi=\frac{2\alpha}{3}$ and $\omega=1$. \\

The main result of this subsection is:   

\begin{theorem}
Let $f:H\to\mathbb{R}$ be $\mu$-strongly convex. Pick $m\ge \frac{1}{2}$, and consider the system \eqref{System: HBF-H}, where
$$\alpha\in \left(0,\, 3\sqrt{\frac{\mu m}{2}} \right],\quad \beta > 0,\quad \gamma = m\left( 1 + \frac{1}{3}\alpha\beta \right) + \frac{1}{3}\alpha\beta.$$
If $v_0 = -\beta\nabla f(x_0)$, we have, for every $t\ge 0$,
$$f\left( x(t) \right) - f^* \le \left( 1 + \frac{6}{\alpha\beta} \right) \left( f(x_0) - f^* \right) e^{-\frac{2}{3}\alpha t}.$$
Moreover,
$$ \int_0^{\infty} e^{\frac{2\alpha}{3}t}\| \nabla f\left(x(t)\right) \|^2 dt < \infty. $$
\end{theorem}

\begin{proof}
Let $W_2$ be the realization of the function $W$, defined by \eqref{E: W}, with $\xi=\frac{2\alpha}{3}$ and $\omega=1$, namely:
$$W_2(t) = \frac{1}{2}\left\| \dot{x} + \beta\nabla f(x) + \frac{2\alpha}{3}(x-x^*) \right\|^2
  - \frac{\alpha^2}{9}\| x - x^* \|^2 + \theta\left( f(x) - f^* \right),$$
where $\theta=\gamma-\frac{1}{3}\alpha\beta>0$. From Lemma \ref{Lem: W_bound}, we obtain
\begin{align*}
\dot{W_2} + \frac{2\alpha}{3} W_2
&\le -\beta\left( \theta - \frac{1}{2} \right)\| \nabla f(x) \|^2 
   + \frac{2\alpha\theta}{3} \left( f(x) - f^* \right)
   - \frac{2\alpha\theta}{3} \langle \nabla f(x), x-x^* \rangle \\
&\quad + \frac{2\alpha^3}{27}\left( 1 + \frac{1}{3}\alpha\beta \right)\| x-x^* \|^2.
\end{align*}
By strong convexity, we have
$$ \langle\nabla f(x),x-x^* \rangle \ge f(x)-f^* + \frac{\mu}{2}\|x-x^*\|^2.$$
Using $\theta = \gamma - \frac{1}{3}\alpha\beta = m\left( 1 + \frac{1}{3}\alpha\beta \right)$ and $\alpha\le 3\sqrt{\frac{\mu m}{2}}$, we obtain
$$ \dot{W_2} + \frac{2\alpha}{3} W_2 + \beta\left( \theta - \frac{1}{2} \right)\| \nabla f(x) \|^2 
\le - \frac{2\alpha}{27}\left( 1 + \frac{1}{3}\alpha\beta \right)\left( \frac{9}{2}\mu m - \alpha^2 \right) \le 0. $$
Since
$$ \beta>0 \qbox{and} \theta = m\left( 1 + \frac{1}{3}\alpha\beta \right) > m \ge \frac{1}{2}, $$
it follows that
$$ W_2(t)\le W_2(0) e^{-\frac{2}{3}\alpha t} \qbox{and} \int_0^{\infty} e^{\frac{2\alpha}{3}t}\| \nabla f\left(x(t)\right) \|^2 dt < \infty. $$
Using $v_0 = -\beta\nabla f(x_0)$ and $\theta= m\left( 1 + \frac{1}{3}\alpha\beta \right)$, we deduce that
\begin{equation}\label{E: V_0_bound_SC_TMM}
\begin{aligned}
W_2(0) 
&= \frac{\alpha^2}{9}\| x_0 - x^* \|^2 + m\left( 1 + \frac{1}{3}\alpha\beta \right)\left( f(x_0) - f^* \right) \\
&\stackrel{(i)}{\le} \frac{\mu m}{2}\| x_0 - x^* \|^2 + m\left( 1 + \frac{1}{3}\alpha\beta \right)\left( f(x_0) - f^* \right) \\
&\stackrel{(ii)}{\le} m\left( 2 + \frac{1}{3}\alpha\beta \right)\left( f(x_0) - f^* \right),
\end{aligned}
\end{equation}
where $(i)$ is due to $\alpha\le 3\sqrt{\frac{\mu m}{2}}$ and $(ii)$ is due to strong convexity. Strong convexity also implies that $f(x) - f^*\ge \frac{\mu}{2}\|x-x^*\|^2$ for every $x\in H$. Remark \ref{Rem: func_bound} shows that 
$$ f\left(x(t) \right) - f^* 
\le \frac{\mu W_2(t)}{\mu\theta - \frac{2}{9}\alpha^2}
\le \frac{W_2(t)}{\theta-m}
= \frac{3W_2(t)}{m \alpha \beta}
\le \frac{3W_2(0)}{m\alpha\beta}e^{-\frac{2}{3}\alpha t}, $$
since $\alpha\le 3\sqrt{\frac{\mu m}{2}}$ and $\theta= m\left( 1 + \frac{1}{3}\alpha\beta \right)$. This, together with \eqref{E: V_0_bound_SC_TMM}, gives the desired result. 
\end{proof}

\begin{remark}
Setting $m=1$ and $\alpha = 3\sqrt{\frac{\mu}{2}}$ so that $\gamma = 1 + \sqrt{2\mu}\beta$, one obtains $f(x) - f^* = \mathcal{O}\left( e^{-\sqrt{2\mu}t} \right)$ for the system
\begin{equation}\label{E: HBF-H-QSC}
\ddot{x} + 3\sqrt{\frac{\mu}{2}} \dot{x} + \beta\nabla^2 f(x) \dot{x} + \left( 1 + \sqrt{2\mu}\beta \right)\nabla f(x) = 0.
\end{equation}
The same rate was derived in \cite[Theorem 3.2]{Aujol_2022} for the heavy-ball system 
\begin{equation}\label{E: HBF-QSC}
\ddot{x} + 3\sqrt{\frac{\mu}{2}} \dot{x} + \nabla f(x) = 0.
\end{equation}
System \eqref{E: HBF-H-QSC} can be  considered as a high-resolution extension of system \eqref{E: HBF-QSC}, and explains the acceleration mechanism behind the algorithm developed in \cite{Aujol_2022}. On the other hand, due to inclusion of a Hessian-driven damping term, discretization of system \eqref{E: HBF-H-QSC} is made much easier and the resultant algorithm is only composed of two steps, in contrast with the four-step scheme for system \eqref{E: HBF-QSC} in \cite{Aujol_2022}.    
\end{remark}

\begin{remark}
Setting $m=2$ and $\alpha = 3\sqrt{\mu}$ so that $\gamma = 2 + 3\sqrt{\mu}\beta$, one gets $f(x) - f^* = \mathcal{O}\left( e^{-2\sqrt{\mu}t} \right)$ for the system
\begin{equation}\label{E: HBF-H-fastest}
\ddot{x} + 3\sqrt{\mu}\dot{x} + \beta\nabla^2 f(x) \dot{x} + \left( 2 + 3\sqrt{\mu}\beta \right)\nabla f(x) = 0.
\end{equation}
A similar system was studied in \cite{Yurtsever_2025} (posted online shortly after the first version of this work) from a different perspective. We shall prove in the next section that this value of $\gamma$ results in a stable discretization that allows us to transpose this factor $2$ gain to the algorithmic setting, explaining the speedup of TMM \cite{Scoy_2018} and ITEM \cite{Taylor_2023_optimal}.  
\end{remark}

\subsection{Convexity, quadratic growth and a unique minimizer} \label{SS:quadunique}

In this subsection, we apply the idea of time rescaling to accelerate the convergence of system \eqref{System: HBF-H} outside of the strongly convex setting. We consider here the case of convex functions with $\mu$-quadratic growth 
$$f(x)-f^*\ge \frac{\mu}{2}\dist(x,S)^2,$$ 
and having a unique minimizer.   

Our analysis relies on the energy function $W$, defined in \eqref{E: W}, but now using $\omega=1$ and $\xi=\frac{2+r}{3+r}\alpha$ for an appropriate value of $r$ which we will select in a moment. We begin by writing
\begin{equation}\label{E: V_QG}
\begin{aligned}
W_3(t) 
&= \frac{1}{2}\left\| \dot{x} + \beta\nabla f(x) + \left(\frac{2+r}{3+r}\right)\alpha(x-x^*) \right\|^2 
 - \frac{(2+r)\alpha^2}{2(3+r)^2}\| x - x^* \|^2 \\
&\quad + \theta\left( f(x) - f^* \right),
\end{aligned}
\end{equation}
where $\theta=\gamma-\frac{\alpha\beta}{3+r}>0$ and $r>0$. Setting $\omega=1$ and $\xi=\frac{2+r}{3+r}\alpha$ in Lemma \ref{Lem: W_bound}, we obtain 
\begin{align*}
\dot{W_3} + \frac{2\alpha }{3+r}\,W_3
&\le -\beta\left( \theta - \frac{1}{2} \right)\| \nabla f(x) \|^2 
  - \frac{(2+r)\alpha\theta}{3+r} \langle \nabla f(x), x-x^* \rangle \\
&\quad  + \underbrace{ \frac{2\alpha\theta}{3+r}\left( f(x) - f^* \right) }_{\mathbf{I}} + \underbrace{ \frac{(2+r)\alpha^3}{2(3+r)^3} \left[ 2(1+r) + \frac{2+r}{3+r}\alpha\beta \right]\| x-x^* \|^2 }_{\mathbf{II}}.
\end{align*}
Since $f$ is convex with quadratic growth and has a unique minimizer, we have 
$$ \langle \nabla f(x), x-x^* \rangle \ge \max\left\{ f(x) - f^*,\, \frac{\mu}{2}\| x-x^* \|^2 \right\}, $$
so that
$$ \mathbf{I} \le \frac{2\alpha\theta}{3+r} \langle \nabla f(x), x-x^* \rangle, $$
and
$$ \mathbf{II} \le \frac{(2+r)\alpha^3}{\mu(3+r)^3} \left[ 2(1+r) + \frac{2+r}{3+r}\alpha\beta \right] \langle \nabla f(x), x-x^* \rangle. $$
Whence, we deduce that
$$\dot{W_3} + \frac{2\alpha}{3+r}W_3 + \beta\left( \theta - \frac{1}{2} \right)\| \nabla f(x) \|^2
\le -\frac{\alpha R}{\mu(3+r)} \langle \nabla f(x), x-x^* \rangle,$$
where
$$ R= r\mu\theta - \frac{\alpha^2(2+r)}{(3+r)^2}\left( 2(1+r) + \frac{2+r}{3+r}\alpha\beta \right). $$
If we guarantee $R\ge 0$, then we will have
\begin{equation} \label{E:W3}
\dot{W_3} + \frac{2\alpha}{3+r}W_3 + \beta\left( \theta - \frac{1}{2} \right)\| \nabla f(x) \|^2 \le 0.
\end{equation} 
Ideally, we would want to maximize $\frac{2\alpha}{3+r}$, subject to $r>0$ and $R\ge 0$. The constraint $R\ge 0$ can be expressed as 
\begin{equation} \label{E:Rge0}
\left( \frac{\alpha}{3+r} \right)^2
\le \frac{ r\mu\theta }{ (2+r) \left( 2(1+r) + \frac{2+r}{3+r}\alpha\beta \right) }.
\end{equation}
To make the rate independent from $\beta$ (in line with our previous results as well as the literature on the subject), we consider the problem
$$\sup_{r>0}\,\inf_{\beta>0}\,\left[ \frac{ r\mu\theta }{ (2+r) \left( 2(1+r) + \frac{2+r}{3+r}\alpha\beta \right)}\right]=\sup_{r>0}\,\left[ \frac{ r\mu\theta }{ 2(1+r)(2+r)}\right],$$
whose solution is $r=\sqrt{2}$. We shall use this value of $r$,\footnote{We acknowledge that this leaves room for improvement in this choice, but $r=\sqrt{2}$ is sufficiently good to obtain our desired improvement on the convergence rate, both in continuous time and in the algorithmic setting discussed in the next section.} and then enforce \eqref{E:Rge0} by limiting the possible values for $\alpha$. More precisely, we have
\begin{align*}
R & = \sqrt{2}\mu\theta - \frac{\alpha^2(2+\sqrt{2})}{(3+\sqrt{2})^2}\left( 2(1+\sqrt{2}) + \frac{2+\sqrt{2}}{3+\sqrt{2}}\alpha\beta \right) \\
& = \sqrt{2}\left[\mu\theta-\alpha^2\frac{(2+\sqrt{2})^2}{(3+\sqrt{2})^2}\left(1+\frac{\alpha\beta}{2+3\sqrt{2}}\right)\right] \\
& = \sqrt{2}\left(1+\frac{\alpha\beta}{2+3\sqrt{2}}\right)\left[\mu m-\alpha^2\left(\frac{2+\sqrt{2}}{3+\sqrt{2}}\right)^2\right] \\
& = \sqrt{2}\left(1+\frac{\alpha\beta}{2+3\sqrt{2}}\right)\left(\frac{2+\sqrt{2}}{3+\sqrt{2}}\right)^2\left[\mu m\left(2-\frac{\sqrt{2}}{2}\right)^2-\alpha^2\right],
\end{align*}
where we have written $m=\theta\left(1+\frac{\alpha\beta}{2+3\sqrt{2}}\right)^{-1}$, which imposes the relationship
$$\gamma = m\left( 1 + \tfrac{\alpha\beta}{2+3\sqrt{2}} \right) + \tfrac{\alpha\beta}{3+\sqrt{2}}. $$
We are now in a position to prove the main result of this subsection, namely:

\begin{theorem}\label{Thm: QG_dynamics}
Let $f:H\to\R$ be convex, have quadratic growth with constant $\mu>0$ and admit a unique minimizer. Pick $m\ge \frac{1}{2}$, and consider the system \eqref{System: HBF-H}, where 
$$ \alpha\in\left( 0, \left( 2 - \tfrac{\sqrt{2}}{2} \right)\sqrt{\mu m} \right],\quad \beta > 0,\quad \gamma = m\left( 1 + \tfrac{\alpha\beta}{2+3\sqrt{2}} \right) + \tfrac{\alpha\beta}{3+\sqrt{2}}. $$
If $v_0 = -\beta\nabla f(x_0)$, we have, for every $t\ge 0$,
$$ f\big( x(t) \big) - f^* \le (1+\sqrt{2}) \big( f(x_0) - f^* \big) e^{-\frac{2}{3+\sqrt{2}}\alpha t }. $$
Moreover,
$$ \int_0^{\infty} e^{\frac{2\alpha}{3+\sqrt{2}}t}\| \nabla f\left(x(t)\right) \|^2 dt < \infty. $$
\end{theorem}

\begin{proof}
Since $\beta>0$ and
$$\theta = \gamma - \tfrac{\alpha\beta}{3+\sqrt{2}} = m\left( 1 + \tfrac{\alpha\beta}{2+3\sqrt{2}} \right) > m \ge \frac{1}{2}, $$
inequality \eqref{E:W3} with $r=\sqrt{2}$ gives
$$ W_3(t) \le W_3(0)e^{-\frac{2}{3+\sqrt{2}}\alpha t} \qbox{and} \int_0^{\infty} e^{\frac{2\alpha}{3+\sqrt{2}}t}\| \nabla f\left(x(t)\right) \|^2 dt < \infty. $$
Using $v_0 = -\beta\nabla f(x_0)$ and $\|x_0-x^*\|^2\le \frac{2}{\mu}(f(x_0-f^*)$ in \eqref{E: V_QG}, we obtain 
\begin{align*}
W_3(0)
&= \frac{(1+r)(2+r)}{2(3+r)^2}\alpha^2 \| x_0 - x^* \|^2
  + \theta\left( f(x_0) - f^* \right) \\
&\le \left[ \frac{(1+r)(2+r)}{(3+r)^2} \left( \frac{\alpha^2}{\mu} \right) + \theta \right]\left( f(x_0) - f^* \right) \\
&\le 
   \left( \frac{1+r}{2+r}m + \theta \right) \left( f(x_0) - f^* \right) \\
&= \left( \frac{m}{\sqrt{2}} + \theta \right)\left( f(x_0) - f^* \right) \quad (\text{since }r=\sqrt{2}) \\
&= \left( \frac{\sqrt{2} + 1 + \frac{\alpha\beta}{3+\sqrt{2}} }{\sqrt{2} + \frac{\alpha\beta}{3+\sqrt{2}} }\right)\theta\left( f(x_0) - f^* \right).
\end{align*}
Since $f$ has quadratic growth and a unique minimizer, it satisfies $f(x) - f^*\ge \frac{\mu}{2}\|x-x^*\|^2$ for every $x\in H$. From Remark \ref{Rem: func_bound}, we deduce that 
$$ f\left(x(t) \right) - f^* 
\le \frac{\mu W_3(t)}{\mu\theta - \frac{2+r}{(3+r)^2}\alpha^2 }
\le \frac{W_3(t)}{\theta-\frac{m}{2+r} }
= \frac{\sqrt{2} + \frac{\alpha\beta}{3+\sqrt{2}} }{ 1 + \frac{\alpha\beta}{3+\sqrt{2}} } \left(\frac{W_3(0)}{\theta}\right) e^{-\frac{2}{3+\sqrt{2}}\alpha t}, $$
and the result follows.
\end{proof}

\begin{remark}
Setting $m=1$ and $\alpha= \left( 2 - \frac{\sqrt{2}}{2} \right)\sqrt{\mu}$ so that $\gamma = 1 + \frac{1}{2}\sqrt{\mu}\beta$, one obtains $f(x) - f^* = \mathcal{O}\left( e^{-(2-\sqrt{2})\sqrt{\mu}t} \right)$ for the system
$$ \ddot{x} + \left( 2 - \frac{\sqrt{2}}{2} \right)\sqrt{\mu}\dot{x} + \beta\nabla^2 f(x) \dot{x} + \left( 1 + \frac{1}{2}\sqrt{\mu}\beta \right) \nabla f(x) = 0.$$ 
The same rate was derived in \cite[Corollary 1]{Aujol_2023} for the heavy-ball system, which corresponds to $\beta = 0$. If we set $m=2$ and $\alpha=(2\sqrt{2}-1)\sqrt{\mu}$ so that $\gamma = 2 + \sqrt{\mu}\beta$. In this case, one gets $f(x) - f^* = \mathcal{O}\left( e^{-2(\sqrt{2}-1)\sqrt{\mu}t} \right)$ for the system
$$\ddot{x} + (2\sqrt{2}-1)\sqrt{\mu} \dot{x} + \beta\nabla^2 f(x) \dot{x} + (2+\sqrt{\mu}\beta) \nabla f(x) = 0.$$   
This rate improves the one in \cite[Corollary 1]{Aujol_2023} by a factor of $\sqrt{2}$, and can be transposed to the algorithmic setting, as we shall see later. 
\end{remark}

\subsection{Polyak-\L ojasiewicz inequality}  \label{SS:quadnonunique}
In this subsection, we derive a linear convergence rate for smooth functions satisfying the Polyak-\L ojasiewicz (P{\L}) inequality with constant $\mu>0$: 
\begin{equation}\label{E: PL}
\|\nabla f(x)\|^2 \ge 2\mu\big(f(x)-f^*\big),\quad\forall x\in H,
\end{equation}
even without requiring convexity. The relationship between the P{\L} inequality and quadratic growth is explained by Corollary \ref{C:PLtoQG} and Proposition \ref{P:QGtoPL}. Consequently, Theorem \ref{T:onlyPL} below can be expressed in terms of quadratic growth, with an appropriate modification of the constants.

As before, we shall use the energy function $W$, defined in \eqref{E: W}, but we now set $\xi=0$, which gives $\theta=\gamma-\alpha\beta$.

We have the following:

\begin{theorem} \label{T:onlyPL}
Let $f:H\to\R$ satisfy the Polyak-\L ojasiewicz inequality \eqref{E: PL} with constant $\mu>0$. Pick $\theta>0$, and consider the system \eqref{System: HBF-H}, where 
$$\beta>0,\quad \gamma=\theta+\alpha\beta. $$
If $v_0 = -\beta\nabla f(x_0)$, we have, for every $t\ge 0$,
$$ f\big( x(t) \big) - f^* \le \big( f(x_0) - f^* \big) e^{-2\delta t },$$
where $\delta:=\min\{\alpha,\mu\beta\}$.
\end{theorem}

\begin{proof}
Let $W_4$ be the realization of the function $W$, defined by \eqref{E: W}, with $\xi=0$, namely:
$$W_4(t) = \frac{1}{2}\left\| \dot{x} + \beta\nabla f(x) \right\|^2 + \theta\left( f(x) - f^* \right),$$
where $\theta=\gamma-\alpha\beta>0$. We shall not use Lemma \ref{Lem: W_bound} this time. Instead, we first compute
\begin{align*}
\dot{W_4} 
&= \left\langle \dot{x} + \beta\nabla f(x), \ddot{x} + \beta\nabla^2 f(x)\dot{x} \right\rangle
  + \theta \langle \nabla f(x), \dot{x} \rangle \\
&= \left\langle \dot{x} + \beta\nabla f(x), -\alpha\dot{x} - \gamma\nabla f(x) \right\rangle
  + \theta \langle \nabla f(x), \dot{x} \rangle \\
&= -\alpha \| \dot{x} + \beta\nabla f(x) \|^2
   - \underbrace{(\gamma-\alpha\beta-\theta)}_{=0}\langle \nabla f(x), \dot{x} \rangle
   - \underbrace{(\gamma-\alpha\beta)}_{=\theta}\beta\| \nabla f(x) \|^2 \\
&= -\alpha \| \dot{x} + \beta\nabla f(x) \|^2 - \beta\theta \| \nabla f(x) \|^2. 
\end{align*}
Using the Polyak-\L ojasiewicz inequality \eqref{E: PL}, we deduce that
$$\dot{W_4}
\le -\alpha \| \dot{x} + \beta\nabla f(x) \|^2 - 2\mu\beta\theta \left( f(x) - f^* \right) 
\le - 2\delta\,W_4.$$
As a result,
$$ W_4(t) \le W_4(0) e^{-2\delta t}. $$
Since $v_0 = -\beta\nabla f(x_0)$, we have $W_4(0)=\theta\left( f(x_0) - f^* \right)$, whence
$$ f\left( x(t) \right) - f^* \le \frac{W_4(t)}{\theta} \le \frac{W_4(0)}{\theta}e^{-2\delta t} \le \left( f(x_0) - f^* \right)e^{-2\delta t}, $$
as claimed.
\end{proof}

\begin{remark}
As shown in Corollary \ref{C:PLtoQG}, the P{\L} inequality implies quadratic growth. Theorem \ref{T:onlyPL} above is different from Theorem \ref{Thm: QG_dynamics}, since it neither requires convexity nor uniqueness of the minimizer. 
\end{remark}

\section{Inertial algorithms}\label{Sec: inertial_algorithms}
In this section, we study accelerated gradient algorithms and establish linear convergence rates for the corresponding function values. These algorithms are obtained by explicit discretization of the inertial system \eqref{System: HBF-H}. In line with the results in Section \ref{Sec: inertial_dynamics}, we recover the speedup upon Nesterov's method by a factor of $\sqrt{2}$ and then $2$ for smooth strongly convex functions. This is given in Subsections \ref{SS:Algorithms_sqrt2} and \ref{SS:Algorithms_2}, respectively. We analyze the case of convex functions with quadratic growth and a unique minimizer in \ref{SS:Algorithms_QG}, establishing a factor $\sqrt{2}$ improvement on the state-of-the-art convergence rate of the algorithm in \cite{Aujol_2023}. Finally, the case of multiple minimizers is discused in \ref{SS:Algorithms_nonunique}.

\subsection{Discretization of \texorpdfstring{\eqref{System: HBF-H}}{} and some useful estimations}
In this subsection, we show the discretization procedures of obtaining the associated inertial algorithms for system \eqref{System: HBF-H}, and present some useful estimations that will be used in the subsequent analysis.

Recall that $f$ is $L$-smooth. We define the step size $h:=\frac{1}{\sqrt{L}}$, set $\beta=h$ and approximate the dynamics \eqref{System: HBF-H} by \eqref{E:Algorithm_family}, which we can rewrite in the more familiar form
\begin{equation}\label{Algo: AGM}\tag{AGM}
\left\{
\begin{array}{rcl}
y_{k+1} & = & x_k - h^2 \nabla f(x_k),\\[2pt]
x_{k+1} & = & y_{k+1} + \frac{1}{1 + \alpha h}( y_{k+1} - y_k ) + \left( \frac{\gamma}{1+\alpha h} - 1 \right) ( y_{k+1} - x_k ).
\end{array}
\right.
\end{equation}
Note that the iterates satisfy
\begin{equation}\label{E: iterates}
(x_{k+1}-y_{k+1}) - (x_k-y_k) = -\alpha h(x_{k+1}-x_k) - \gamma h^2 \nabla f(x_k).
\end{equation}

\begin{remark}\label{Rem: x_0}
The algorithm \eqref{Algo: AGM} can be initialized with given $x_0$ and $y_0$. From \eqref{E: iterates}, one has  
$$ (1+\alpha h)(x_1-x_0) = x_0 - y_0 - (\gamma+1)h^2\nabla f(x_0), $$
which implies that we can determine $y_0$ by specifying $x_0$ and $x_1$:
$$ y_0 = -(1+\alpha h)(x_1-x_0) - (\gamma+1)h^2\nabla f(x_0) + x_0. $$
This observation will be useful in the subsequent convergence analysis. 
\end{remark}

We have the following:
\begin{lemma}\label{Lem: func_bound}
Let $f$ be $\mu$-strongly convex and $L$-smooth. Let $(x_k)$ and $(y_k)$ be generated by \eqref{Algo: AGM}. Define
\begin{equation}\label{E: psi_k}
\psi_k := f(x_k) - f^* - (h^2/2)\| \nabla f(x_k) \|^2.
\end{equation}
Then, we have
$$ \langle \nabla f(x_k), x_{k+1}-x_k \rangle 
\ge \frac{1}{1+\alpha h}\left[ \psi_k - \psi_{k-1} - \gamma h^2 \| \nabla f(x_k) \|^2 + \frac{\mu}{2(1-\mu h^2)}\| y_{k+1} - y_k \|^2 \right], $$
and
$$ \langle \nabla f(x_k), x_k - x^* \rangle 
\ge \psi_k + \frac{h^2}{2} \| \nabla f(x_k) \|^2 + \frac{\mu}{2}\| x_k-x^* \|^2. $$
\end{lemma}

\begin{proof}
We begin by using \eqref{Algo: AGM}, to obtain 
$$ x_k - x_{k-1} = (1+\alpha h)(x_{k+1}-x_k) + h^2\left( \nabla f(x_k) - \nabla f(x_{k-1}) \right) + \gamma h^2 \nabla f(x_k). $$
Since $f$ is $\mu$-strongly convex and $L$-smooth, we have
\begin{align}
\nonumber
&\quad f(x_k) - f(x_{k-1}) \\
\label{E: interpolation_condition}
&\le \langle \nabla f(x_k), x_k - x_{k-1} \rangle - \frac{h^2}{2}\| \nabla f(x_k) - \nabla f(x_{k-1}) \|^2 - \frac{\mu\| y_{k+1} - y_k \|^2}{2(1-\mu h^2)}\\
\nonumber
&= (1+\alpha h)\langle \nabla f(x_k), x_{k+1}-x_k \rangle
  + h^2\left\langle \nabla f(x_k), \nabla f(x_k) - \nabla f(x_{k-1}) \right\rangle\\
\nonumber
&\quad + \gamma h^2 \| \nabla f(x_k) \|^2 
  - \frac{h^2}{2}\| \nabla f(x_k) - \nabla f(x_{k-1}) \|^2
  - \frac{\mu}{2(1-\mu h^2)}\| y_{k+1} - y_k \|^2 \\
\nonumber
&= (1+\alpha h)\langle \nabla f(x_k), x_{k+1}-x_k \rangle 
  + \frac{h^2}{2}\| \nabla f(x_k) \|^2
  - \frac{h^2}{2}\| \nabla f(x_{k-1}) \|^2 \\
\nonumber
&\quad + \gamma h^2 \| \nabla f(x_k) \|^2  
  - \frac{\mu}{2(1-\mu h^2)}\| y_{k+1} - y_k \|^2,  
\end{align}
so that
\begin{align*}
\psi_k - \psi_{k-1} 
&\le (1+\alpha h)\langle \nabla f(x_k), x_{k+1}-x_k \rangle + \gamma h^2 \| \nabla f(x_k) \|^2 
 - \frac{\mu}{2(1-\mu h^2)}\| y_{k+1} - y_k \|^2.
\end{align*}
This gives the first claim. Likewise,
$$ \langle \nabla f(x_k), x_k - x^* \rangle 
\ge f(x_k) - f^* + \frac{\mu}{2}\| x_k-x^* \|^2 
= \psi_k + \frac{h^2}{2} \| \nabla f(x_k) \|^2 + \frac{\mu}{2}\| x_k-x^* \|^2, $$
which allows us to conclude.
\end{proof}

\begin{remark}
In case $f$ is only convex, the results still hold with $\mu=0$.
\end{remark}

\begin{remark}
Note that $\psi_k \ge f(y_{k+1})-f^*\ge 0$ due to the convexity and $L$-smoothness of $f$.
\end{remark}

\begin{remark}
The inequality \eqref{E: interpolation_condition} holds for smooth strongly convex functions, and is referred to as the interpolation condition in the context of Performance Estimation Problems \cite[Theorem 4]{Talor_2017_SC}.
\end{remark}

The following inequalities will facilitate the forthcoming analysis.

\begin{lemma}
Let $a,b,c>0$. For every $x>0$, we have
\begin{equation}\label{E: ineqn_frac}
\frac{ (1+ax)(1+bx) }{1+cx} \le 1 + x\max\left\{ a+b-c, \frac{ab}{c} \right\},
\end{equation}
and
\begin{equation}\label{E: ineqn_frac_quad}
\frac{1+ax+bx^2}{1+cx} \le 1 + x\max\left\{ a-c, \frac{b}{c} \right\}.
\end{equation}
\end{lemma}

\begin{proof}
Let $u = \max\left\{ a+b-c, \frac{ab}{c} \right\}$. Then,
$$ (1+cx)(1+ux) = 1 + (u+c)x + ucx^2 \ge 1 + (a+b)x + abx^2 = (1+ax)(1+bx), $$
so that \eqref{E: ineqn_frac} holds. Likewise, we can also prove \eqref{E: ineqn_frac_quad}.
\end{proof}

\subsection{Strong convexity: factor \texorpdfstring{$\sqrt{2}$}{} acceleration} \label{SS:Algorithms_sqrt2}
In this subsection, we present the convergence results which improve from that of Nesterov's method \eqref{E: NAG-SC} to those of \cite{Luo_2021,Park_2023}, with a factor of $\sqrt{2}$. As one will see, the improvement can be attributed to an increasing value of $\gamma$ in \eqref{Algo: AGM}. 

Our analysis centers around the energy sequence $(E_k)_{k\ge 1}$, defined by
\begin{equation}\label{E: E_k_SC}
E_k = \frac{1}{2}\| \phi_k \|^2 + \theta h^2 \psi_{k-1},
\end{equation}
where $\psi_k$ is given by \eqref{E: psi_k}, and
$$ \phi_k = (x_k-y_k) + \frac{ \frac{1}{2}\alpha h }{ 1 + \frac{1}{2}\alpha h }(x_k-x^*),
\quad
\theta = \frac{ \gamma-\frac{1}{2}\alpha h }{1 + \frac{1}{2}\alpha h}. $$

\begin{remark}
The $(E_k)_{k\ge 1}$ sequence was introduced in \cite{Park_2023}. In particular, the term $\psi_{k-1}$ will facilitate the Lyapunov analysis. It provides an upper bound for the function values, namely: $f(y_k)-f^*\le\psi_{k-1}$.
\end{remark}

We have the following lemma, whose proof we postpone to Appendix \ref{Proof: Lem: E_k_bound_SC}:

\begin{lemma}\label{Lem: E_k_bound_SC}
Let $f$ be strongly convex. Let $(x_k)$ and $(y_k)$ be generated by \eqref{Algo: AGM} with $\gamma\in[1+\alpha h, 2+\alpha h]$. Consider $(E_k)_{k\ge 1}$ defined by \eqref{E: E_k_SC}. Then, 
$$ \left( 1 + \frac{1}{2}\alpha h \right)E_{k+1} - E_k + \frac{1}{4}\alpha h^5 \| \nabla f(x_k) \|^2 
\le \frac{\frac{1}{16}\alpha h^3 (\alpha^2-4\mu\gamma)}{\left( 1+\frac{1}{2}\alpha h \right)^2}\| x_k - x^* \|^2.  $$
\end{lemma}

We are now in a position to present the main result of this subsection:

\begin{theorem}\label{Thm: algo_SC_OGM}
Let $f:H\to\R$ be $\mu$-strongly convex and $L$-smooth. Pick $m\in[1,2]$. Let $(x_k)$ and $(y_k)$ be generated by \eqref{Algo: AGM}, where
$$ h=1/\sqrt{L},\qbox{}\alpha=2\sqrt{\mu m}, \qbox{and} \gamma=m+\alpha h. $$
If $y_0 = -\left(\gamma - \frac{1}{2}\alpha h \right)h^2 \nabla f(x_0)+x_0$, then for every $k\ge 1$,
$$ f(y_k) - f^* \le \left[ \frac{\mu}{2}\| x_0 - x^* \|^2 + \left( f(x_0) - f^* \right) \right] (1 + \rho)^{-k+1}, $$
where $\rho = \sqrt{qm}$ and $q=\mu/L$. Moreover, $ \sum_{k=0}^\infty (1+\rho)^k \| \nabla f(x_k) \|^2 < \infty $. In particular, $\| \nabla f(x_k) \|^2=o\big((1+\rho)^{-k}\big)$.
\end{theorem}

\begin{proof}
Note that $\alpha^2=4\mu m \le 4\mu\gamma$. It follows from Lemma \ref{Lem: E_k_bound_SC} that
$$ (1 + \rho )E_{k+1} - E_k + \frac{1}{2}\rho h^4 \| \nabla f(x_k) \|^2 \le 0,\qbox{with} \rho = \sqrt{qm}. $$
Hence,
$$ E_k \le (1 + \rho)^{-k+1}E_1 \qbox{and} \sum_{k=0}^\infty (1+\rho)^k \| \nabla f(x_k) \|^2 < \infty. $$
Recalling Remark \ref{Rem: x_0}, we can verify that the choice of $y_0 = - \left(\gamma - \frac{1}{2}\alpha h \right)h^2 \nabla f(x_0) + x_0$ guarantees
\begin{align*}
\phi_1 
&= (x_1-y_1) + \frac{ \frac{1}{2}\alpha h }{ 1 + \frac{1}{2}\alpha h }(x_1-x^*) \\
&= \underbrace{ \left( \frac{1+\alpha h}{ 1 + \frac{1}{2}\alpha h } \right)(x_1-x_0) + h^2\nabla f(x_0) }_{=0}
  + \frac{ \frac{1}{2}\alpha h }{ 1 + \frac{1}{2}\alpha h }(x_0-x^*) \\
&= \frac{ \frac{1}{2}\alpha h }{ 1 + \frac{1}{2}\alpha h }(x_0-x^*).
\end{align*} 
Hence, in view of \eqref{E: E_k_SC}, we obtain 
\begin{align*}
E_1 
&= \frac{ \frac{1}{8}\alpha^2 h^2 }{ \left( 1 + \frac{1}{2}\alpha h \right)^2 }\| x_0 - x^* \|^2 + \theta h^2 \left( f(x_0) - f^* \right)\\
&\le \frac{ \frac{1}{2}\mu m h^2 }{1 + \frac{1}{2}\alpha h} \| x_0 - x^* \|^2 + \theta h^2 \left( f(x_0) - f^* \right) \\
&\le \theta h^2 \left[ \frac{\mu}{2}\| x_0 - x^* \|^2 + \left( f(x_0) - f^* \right) \right], 
\end{align*}
due to strong convexity and $\frac{m}{1 + \frac{1}{2}\alpha h} \le \theta$. Whence, we deduce that
$$ f(y_k) - f^* \le \psi_{k-1} 
\le \frac{E_k}{\theta h^2}
\le \left[ \frac{\mu}{2}\| x_0 - x^* \|^2 + \left( f(x_0) - f^* \right) \right] (1 + \rho)^{-k+1}, $$
for every $k\ge 1$, as claimed.
\end{proof}

\begin{remark}
Setting $m=1$ so that $\alpha=2\sqrt{\mu}$ and $\gamma=1+\alpha h$, one gets a convergence rate $f(y_k)-f^*=\mathcal{O}\left( (1+\sqrt{q})^{-k} \right)$ for Nesterov-type method \cite{Nesterov_2004}:
\begin{equation*}
\left\{
\begin{array}{rcl}
y_{k+1} & = & x_k - h^2 \nabla f(x_k),\\[2pt]
x_{k+1} & = & y_{k+1} + \frac{1}{1 + \alpha h}( y_{k+1} - y_k ).
\end{array}
\right.
\end{equation*} 
By increasing the value of $\gamma$, one gets a faster convergence rate. In particular, setting $m=2$ so that $\alpha=2\sqrt{2\mu}$ and $\gamma=2+\alpha h$, one has $f(y_k)-f^*=\mathcal{O}\left( (1+\sqrt{2q})^{-k} \right)$ for OGM-SC-type method \cite{Park_2023}: 
\begin{equation*}
\left\{
\begin{array}{rcl}
y_{k+1} & = & x_k - h^2 \nabla f(x_k),\\[2pt]
x_{k+1} & = & y_{k+1} + \frac{1}{1 + \alpha h}( y_{k+1} - y_k ) + \frac{1}{1 + \alpha h} ( y_{k+1} - x_k ).
\end{array}
\right.
\end{equation*} 
The $\sqrt{2}$ improvement in the convergence rate is in line with the results in the continuous time. 
\end{remark}

\begin{remark}
Note that \eqref{Algo: AGM} is obtained using {\it explicit discretization} of the time-rescaled inertial system \eqref{System: HBF-H}. In spite of the promising convergence rates in the continuous time using time rescaling, there exist constraints in recovering these results in the discrete time: the value of $\gamma$ cannot be arbitrarily large and satisfies $\gamma_{\max}\sim 2$ (see also Theorem \ref{T:Algorithm_2} below). The constraints on $\gamma$ is in general due to the $L$-smoothness of $f$ \cite{Attouch_2019_proximal}. 
\end{remark}

\subsection{Strong convexity: factor \texorpdfstring{$2$}{} acceleration}\label{SS:Algorithms_2}
In this subsection, we derive faster convergence rates which range from a factor of $\sqrt{2}$ to $2$, matching the state-of-the-art results. The improved convergence rates are obtained using time rescaling and a new energy sequence.

Our analysis centers around the energy sequence $(E_k)_{k\ge 1}$, defined by
\begin{equation}\label{E: E_k_SC_TMM}
E_k = \frac{1}{2}\| \phi_k \|^2 - \frac{\eta}{2}\| y_k-x^* \|^2 + \theta h^2 \psi_{k-1},
\end{equation}
where $\psi_k$ is given by \eqref{E: psi_k}, and
\begin{align*}
\phi_k &= (x_k-y_k) + \frac{2}{3}\alpha h(x_k-x^*),
\quad \eta = \frac{\frac{2}{9}\alpha^2 h^2}{1 + \alpha h },\\
\theta &= \gamma  - \frac{1}{3}\alpha h\left( \frac{1 + \frac{2}{3}\alpha h}{1+\alpha h} \right).
\end{align*}

\begin{remark}\label{Rem: func_estimate_QG}
Let $f$ be $L$-smooth, have $\mu$-quadratic growth and admit a unique minimizer. Consider an energy sequence $(E_k)_{k\ge 1}$ of the form \eqref{E: E_k_SC_TMM}. Observe that
\begin{align*}
E_k &\ge - \frac{\eta}{2}\left\| y_k - x^* \right\|^2 + \theta h^2 \left[ f(x_{k-1}) - f^* - \frac{h^2}{2}\| \nabla f(x_{k-1}) \|^2 \right] \\
&\ge - \frac{\eta}{2}\left\| y_k - x^* \right\|^2 + \theta h^2 \left( f(y_k) - f^* \right)\\
&\ge \left( \theta h^2 - \frac{\eta}{\mu} \right) \left( f(y_k) - f^* \right).
\end{align*}
If there exists a constant $\rho > 0$ such that
\begin{equation}\nonumber
(1+\rho) E_{k+1} - E_k \le 0,
\end{equation}
and $\theta h^2 > \frac{\eta}{\mu}$, we can deduce that
\begin{equation}\nonumber
f(y_k) - f^* \le \frac{E_k}{\theta h^2 - \frac{\eta}{\mu}} \le \frac{E_1 (1+\rho)^{-k+1}}{\theta h^2 - \frac{\eta}{\mu} },\quad \forall k\ge 1,
\end{equation}
which implies a linear convergence rate for $f(y_k) - f^*$.
\end{remark}

We have the following result, whose proof can be found in Appendix \ref{Proof: Lem: E_k_bound_SC_TMM}:

\begin{lemma}\label{Lem: E_k_bound_SC_TMM}
Let $f$ be strongly convex. Let $(x_k)$ and $(y_k)$ be generated by \eqref{Algo: AGM} with $\gamma\in[1, 2]$ and $\alpha h\in(0,3]$. Consider $(E_k)_{k\ge 1}$ defined by \eqref{E: E_k_SC_TMM}. Then, 
\begin{align*}
\left( 1 + \frac{\frac{2}{3}\alpha h}{1 + \frac{2}{3}\alpha h} \right) E_{k+1} - E_k
  + \frac{ \frac{1}{3}\alpha h^5 \| \nabla f(x_k) \|^2 }{ 1 + \frac{2}{3}\alpha h }
&\le \frac{\frac{2}{27}\alpha h^3 (1+2\alpha h) \left(\alpha^2 - \frac{9}{2}\mu\gamma \right)}{\left( 1+\frac{2}{3}\alpha h \right)(1+\alpha h)}\| x_k - x^* \|^2 \\[2pt]
&\quad + \frac{\frac{1}{9} h^2 \left( \alpha^2 - \frac{9}{2}\mu\gamma \right)}{1 + \alpha h } \| y_{k+1} - y_k \|^2. 
\end{align*}
\end{lemma}

Now we present the main result of this subsection in what follows.

\begin{theorem} \label{T:Algorithm_2}
Let $f:H\to\R$ be $\mu$-strongly convex and $L$-smooth. Let $(x_k)$ and $(y_k)$ be generated by \eqref{Algo: AGM}, where
$$ h=1/\sqrt{L},\qbox{} \gamma\in[1,2],\qbox{and}\alpha=3\sqrt{\mu\gamma/2}. $$
If $y_0 = -\left( \gamma - \frac{ \frac{1}{3}\alpha h }{1+ \frac{2}{3}\alpha h } \right)h^2\nabla f(x_0) + x_0 $, then for every $k\ge 1$,
$$ f(y_k) - f^* \le C (1+\rho)^{-k+1}, $$
where $\rho = \frac{\sqrt{2q\gamma}}{1+\sqrt{2q\gamma}}$, $q=\mu/L$, and
$$ C = \frac{1}{2-\sqrt{2}} \left( 3 + \sqrt{2/(q\gamma)} \right) \left[ \mu \| x_0 - x^* \|^2 + \left( f(x_0) - f^* \right)  \right]. $$
Moreover, $\sum_{k=0}^\infty (1+\rho)^k \| \nabla f(x_k) \|^2 < \infty$, and so $\| \nabla f(x_k) \|^2=o\big((1+\rho)^{-k}\big)$. 
\end{theorem}

\begin{proof}
Note that $\alpha h = 3\sqrt{q\gamma/2}\le 3$ since $q\le 1$ and $\gamma\le 2$. It follows from Lemma \ref{Lem: E_k_bound_SC_TMM} that
$$ (1 + \rho )E_{k+1} - E_k + \frac{1}{2}\rho h^4 \| \nabla f(x_k) \|^2 \le 0,\qbox{with} \rho = \frac{\sqrt{2q\gamma}}{1+\sqrt{2q\gamma}}. $$
Hence, 
$$ E_k \le (1 + \rho)^{-k+1}E_1 \qbox{and} \sum_{k=0}^\infty (1+\rho)^k \| \nabla f(x_k) \|^2 < \infty. $$
Recalling Remark \ref{Rem: x_0}, it is easy to verify that the choice of $y_0 = -\left( \gamma - \frac{ \frac{1}{3}\alpha h }{1+ \frac{2}{3}\alpha h } \right)h^2\nabla f(x_0) + x_0 $ leads to
\begin{align*}
\phi_1
&= (x_1-y_1) + \frac{2}{3}\alpha h(x_1-x^*) \\
&= \underbrace{ \left( 1 + \frac{2}{3}\alpha h \right)(x_1-x_0) + h^2 \nabla f(x_0) }_{=0}
  + \frac{2}{3}\alpha h (x_0-x^*)\\
&= \frac{2}{3}\alpha h (x_0-x^*).
\end{align*}
In view of \eqref{E: E_k_SC_TMM}, we deduce that  
\begin{align*}
E_1
&\le \frac{2}{9}\alpha^2 h^2 \| x_0 - x^* \|^2 + \gamma h^2 \left( f(x_0) - f^* \right) \\
&= \gamma h^2 \left[ \mu \| x_0 - x^* \|^2 + \left( f(x_0) - f^* \right)  \right]. 
\end{align*}
On the other hand, according to Remark \ref{Rem: func_estimate_QG}, we start with
\begin{align*}
\theta h^2 - \frac{\eta}{\mu}
&= h^2 \left[ \gamma  - \frac{1}{3}\alpha h\left( \frac{1 + \frac{2}{3}\alpha h}{1+\alpha h} \right) -  \frac{\frac{2\alpha^2}{9\mu}}{1 + \alpha h } \right] \\
&= \frac{\alpha h^3}{1+\alpha h} \left[ \gamma - \frac{1}{3}\left( 1 + \frac{2}{3}\alpha h \right) \right]\quad (\text{since } \alpha^2 = 9\mu\gamma/2 ) \\
&\ge \frac{\alpha h^3}{3(1+\alpha h)}\left[ (3-\sqrt{2})\gamma-1 \right] \quad (\text{since }\alpha h \le 3\sqrt{\gamma/2} )\\
&\ge \frac{(2-\sqrt{2})\gamma\alpha h^3}{3(1+\alpha h)} \quad (\text{since }\gamma\ge 1 ),
\end{align*}
to arrive at
\begin{align*}
f(y_k) - f^*  
&\le \frac{E_1 (1+\rho)^{-k+1}}{\theta h^2 - \frac{\eta}{\mu} } \\ 
&\le \frac{3}{2-\sqrt{2}} \left( 1 + \frac{1}{\alpha h} \right) \left[ \mu \| x_0 - x^* \|^2 + \left( f(x_0) - f^* \right)  \right] (1+\rho)^{-k+1} \\
&= \frac{1}{2-\sqrt{2}} \left( 3 + \frac{\sqrt{2}}{\sqrt{q\gamma}} \right) \left[ \mu \| x_0 - x^* \|^2 + \left( f(x_0) - f^* \right)  \right] (1+\rho)^{-k+1}, 
\end{align*}
as claimed.
\end{proof}

\begin{remark}
For $q$ small enough, one has $\rho\sim \sqrt{2q\gamma}$. Setting $\gamma=1$, one obtains a convergence rate $f(y_k)-f^*=\mathcal{O}\left( (1+\sqrt{2q})^{-k} \right)$, matching that of \cite{Aujol_2022}. Setting $\gamma=2$, the convergence rate can be accelerated by a factor of $\sqrt{2}$ and reaches $f(y_k)-f^*=\mathcal{O}\left( (1+2\sqrt{q})^{-k} \right)$, matching asymptotically the state-of-the-art rates for TMM \cite{Scoy_2018} and ITEM \cite{Taylor_2023_optimal}.    
\end{remark}

\subsection{Convexity, quadratic growth and a unique minimizer}\label{SS:Algorithms_QG}
In this subsection, we show that the time rescaling also works for accelerating the convergence rate for convex functions with quadratic growth and a unique minimizer. In particular, we improve the best rate $\mathcal{O}\left( (1+(2-\sqrt{2})\sqrt{q})^{-k} \right)$, obtained in \cite{Aujol_2023}, by a factor of $\sqrt{2}$.

Our analysis centers around the energy sequence $(E_k)_{k\ge 1}$, defined by
\begin{equation}\label{E: E_k_QG}
E_k = \frac{1}{2}\| \phi_k \|^2 - \frac{\eta}{2}\| y_k-x^* \|^2 + \theta h^2 \psi_{k-1},
\end{equation}
where $\psi_k$ is given by \eqref{E: psi_k}, and
\begin{align*}
\phi_k &= (x_k-y_k) + \xi h(x_k-x^*),
\quad \eta = \frac{\xi(\alpha-\xi) h^2}{(1+\alpha h)(1+2\xi h)},\\
\theta &= \gamma - \frac{(\alpha-\xi)h \left( 1 + \xi h (1+\gamma) \right) }{(1+\alpha h)( 1+2\xi h )},\quad \xi = \frac{2+\sqrt{2}}{3+\sqrt{2}}\alpha.
\end{align*}
The value of $\xi$ comes from the analysis in the continuous time.

We shall use the following result, whose proof can be found in Appendix \ref{Proof: Lem: E_k_bound_QG}:

\begin{lemma}\label{Lem: E_k_bound_QG}
Let $f$ be convex, have $\mu$-quadratic growth and admit a unique minimizer. Let $(x_k)$ and $(y_k)$ be generated by \eqref{Algo: AGM} with $\gamma\in[1, 2]$. Consider $(E_k)_{k\ge 1}$ defined by \eqref{E: E_k_QG}. Then, 
\begin{align*}
&\quad \left( 1 + \frac{ 2(\alpha-\xi) h }{1+2\xi h} \right)E_{k+1} - E_k
 + \frac{ (\alpha-\xi) h^5 \| \nabla f(x_k) \|^2 }{ 1 + 2\xi h } \\
&\le \frac{ (\alpha-\xi) h^3 \left( 1 + \frac{\xi(3\alpha-2\xi)h}{2\xi-\alpha} \right) \left[ \xi(2\xi-\alpha) - \frac{\mu\gamma}{\sqrt{2}} \right]}{(1+\alpha h)\left(1+2\xi h\right)}\| x_k - x^* \|^2,
\end{align*}
where $\xi = \frac{2+\sqrt{2}}{3+\sqrt{2}}\alpha$. 
\end{lemma}

Now we present the main result of this subsection:

\begin{theorem}\label{Thm: algo_QG_OGM}
Let $f:H\to\R$ be convex and $L$-smooth, have $\mu$-quadratic growth and admit a unique minimizer. Let $(x_k)$ and $(y_k)$ be generated by \eqref{Algo: AGM}, where
$$ h=1/\sqrt{L},\qbox{} \gamma\in[1,2],\qbox{and}\alpha=\frac{3+\sqrt{2}}{2+\sqrt{2}}\sqrt{\mu\gamma}. $$
If $y_0 = -\left( \gamma - \frac{\sqrt{q\gamma}}{ (2+\sqrt{2})(1+\sqrt{q\gamma}) } \right)h^2\nabla f(x_0) + x_0$, then for every $k\ge 1$,
$$ f(y_k) - f^* \le \left[ \frac{\sqrt{2}\mu}{2} \| x_0 - x^* \|^2 + \sqrt{2}\left( f(x_0) - f^* \right)  \right] (1+\rho)^{-k+1}, $$
where $\rho = \frac{ (2-\sqrt{2})\sqrt{q\gamma} }{ 1 + 2\sqrt{q\gamma} }$ and $q=\mu/L$. Moreover, $\sum_{k=0}^\infty (1+\rho)^k \| \nabla f(x_k) \|^2 < \infty$.
\end{theorem}

\begin{proof}
Using $\xi = \frac{2+\sqrt{2}}{3+\sqrt{2}}\alpha$, we have
$$ \xi(2\xi-\alpha) - \frac{\mu\gamma}{\sqrt{2}} = \frac{(2+\sqrt{2})(1+\sqrt{2})}{(3+\sqrt{2})^2}\alpha^2 - \frac{\mu\gamma}{\sqrt{2}} = \frac{1}{\sqrt{2}}\left[ \left( \frac{2+\sqrt{2}}{3+\sqrt{2}} \right)^2 \alpha^2 - \mu\gamma \right] = 0, $$
in view of $\alpha=\frac{3+\sqrt{2}}{2+\sqrt{2}}\sqrt{\mu\gamma}$. It follows from Lemma \ref{Lem: E_k_bound_QG} that
$$ (1+\rho)E_{k+1} - E_k + \frac{1}{2}\rho h^4 \| \nabla f(x_k) \|^2 \le 0,\qbox{with} \rho = \frac{ (2-\sqrt{2})\sqrt{q\gamma} }{ 1 + 2\sqrt{q\gamma} }. $$
Hence,
$$ E_k \le (1+\rho)^{-k+1}E_1 \qbox{and} \sum_{k=0}^\infty (1+\rho)^k \| \nabla f(x_k) \|^2 < \infty. $$
Recalling Remark \ref{Rem: x_0}, it is easy to verify that $y_0 = -\left( \gamma - \frac{\sqrt{q\gamma}}{ (2+\sqrt{2})(1+\sqrt{q\gamma}) } \right)h^2\nabla f(x_0) + x_0$ gives 
\begin{align*}
\phi_1 
&= (x_1-y_1) + \frac{2+\sqrt{2}}{3+\sqrt{2}}\alpha h (x_1-x^*) \\
&= (x_1-y_1) + \sqrt{q\gamma}(x_1-x^*) \\
&= \underbrace{ ( 1 + \sqrt{q\gamma} ) (x_1-x_0) + h^2\nabla f(x_0) }_{=0} 
   + \sqrt{q\gamma} (x_0-x^*) \\
&= \sqrt{q\gamma} (x_0-x^*).
\end{align*}
In view of \eqref{E: E_k_QG}, we obtain
\begin{align*}
E_1
&\le \frac{q\gamma}{2} \| x_0 - x^* \|^2 + \gamma h^2 \left( f(x_0) - f^* \right) \\
&= \gamma h^2 \left[ \frac{\mu}{2} \| x_0 - x^* \|^2 + \left( f(x_0) - f^* \right)  \right]. 
\end{align*}
On the other hand, according to Remark \ref{Rem: func_estimate_QG}, and keeping in mind $\alpha^2 = \left( \frac{3+\sqrt{2}}{2+\sqrt{2}} \right)^2 \mu\gamma$, we obtain
\begin{align*}
\theta h^2 - \frac{\eta}{\mu}
&= h^2 \left[ \gamma - \frac{(\alpha-\xi)h \left( 1 + \xi h (1+\gamma) \right) }{(1+\alpha h)( 1+2\xi h )} - \frac{ \xi(\alpha-\xi)/\mu }{(1+\alpha h)(1+2\xi h)} \right] \\
&= h^2 \left[ \gamma - \frac{(\alpha-\xi)h \left( 1 + \xi h (1+\gamma) \right) }{(1+\alpha h)( 1+2\xi h )} - \frac{ \frac{ 2+\sqrt{2} }{(3+\sqrt{2})^2} \frac{\alpha^2}{\mu} }{(1+\alpha h)(1+2\xi h)} \right] \\
&= \frac{\sqrt{2}}{2}\gamma h^2 + \frac{ \frac{\gamma}{2+\sqrt{2}} \left( (\alpha+2\xi)h + 2\alpha\xi h^2 \right) - (\alpha-\xi)h \left( 1 + \xi h (1+\gamma) \right) }{(1+\alpha h)(1+2\xi h)}\\
&= \frac{\sqrt{2}}{2}\gamma h^2 + \frac{ \frac{\gamma (\alpha+2\xi)h}{2+\sqrt{2}} - (\alpha-\xi)h + \xi h^2 \left( \frac{2\alpha\gamma}{2+\sqrt{2}} - (\alpha-\xi)(1+\gamma)  \right) }{(1+\alpha h)(1+2\xi h)} \\
&\ge \frac{\sqrt{2}}{2}\gamma h^2, 
\end{align*}
so that
$$ f(y_k) - f^*  
\le \frac{E_1 (1+\rho)^{-k+1}}{\theta h^2 - \frac{\eta}{\mu} } 
\le \left[ \frac{\sqrt{2}\mu}{2} \| x_0 - x^* \|^2 + \sqrt{2}\left( f(x_0) - f^* \right)  \right] (1+\rho)^{-k+1}, $$
which allows us to conclude.
\end{proof}

\begin{remark}
For $q$ small enough, one has $\rho\sim (2-\sqrt{2})\sqrt{q\gamma}$. Setting $\gamma=1$, one gets a convergence rate $f(y_k)-f^*=\mathcal{O}\left( 1 + (2-\sqrt{2})\sqrt{q} \right)^{-k}$, which matches the one obtained in \cite{Aujol_2023}. Choosing $\gamma=2$, one obtains a faster convergence rate $f(y_k)-f^*=\mathcal{O}\left( 1 + 2(\sqrt{2}-1)\sqrt{q} \right)^{-k}$, which improves on \cite{Aujol_2023} by a factor $\sqrt{2}$.   
\end{remark}

\subsection{P{\L} convex functions with possibly multiple minimizers}\label{SS:Algorithms_nonunique}

In this subsection, we show that the time rescaling can be useful in deriving a faster convergence rate for smooth convex functions that satisfy the Polyak-\L ojasiewicz inequality \eqref{E: PL}.

Our analysis is based on the energy sequence $(E_k)_{k\ge 1}$, defined by
\begin{equation}\label{E: E_k_PL}
E_k = \frac{1}{2}\| \phi_k \|^2 + \gamma h^2 \psi_{k-1},
\end{equation}
where $\phi_k = x_k-y_k$ and $\psi_k$ is given by \eqref{E: psi_k}.

\begin{lemma}\label{Lem: E_k_bound_PL}
Let $f$ be convex. Let $(x_k)$ and $(y_k)$ be generated by \eqref{Algo: AGM}. Consider $(E_k)_{k\ge 1}$ defined by \eqref{E: E_k_PL}. Then, 
$$ (1+\alpha h)E_{k+1}-E_k 
\le - \left[ \gamma(2-\gamma) + (\gamma-1) \alpha h \right] \frac{h^4}{2} \| \nabla f(x_k) \|^2
    + \gamma\alpha h^3 \left( f(x_k) - f^* \right). $$
\end{lemma}

\begin{proof}
By definition of $\phi_k$, we have
$$ \phi_{k+1} = x_{k+1}-y_{k+1} = x_{k+1}-x_k + h^2 \nabla f(x_k), $$
and in view of \eqref{E: iterates},
\begin{align*}
\phi_{k+1} - \phi_k
&= -\alpha h ( x_{k+1} - x_k ) - \gamma h^2 \nabla f(x_k) \\
&= - \frac{\alpha h}{2}\phi_{k+1} 
   - \frac{\alpha h}{2}( x_{k+1} - x_k )
   - \left(\gamma - \frac{1}{2}\alpha h \right) h^2 \nabla f(x_k).
\end{align*}
It follows that
$$ \| \phi_{k+1} - \phi_k \|^2
= \alpha^2 h^2 \| x_{k+1} - x_k \|^2
  + \gamma^2 h^4 \| \nabla f(x_k) \|^2
  + 2\gamma\alpha h^3 \langle \nabla f(x_k), x_{k+1}-x_k \rangle, $$
and
\begin{align*}
\langle \phi_{k+1}-\phi_k, \phi_{k+1} \rangle
&= -\frac{\alpha h}{2}\| \phi_{k+1} \|^2
 - \frac{\alpha h}{2}\| x_{k+1}-x_k \|^2
 - \gamma h^2 \langle \nabla f(x_k), x_{k+1}-x_k \rangle \\
&\quad - \left( \gamma - \frac{1}{2}\alpha h \right) h^4 \| \nabla f(x_k) \|^2.
\end{align*}
Since $\frac{1}{2}\| \phi_{k+1} \|^2 - \frac{1}{2}\| \phi_k \|^2 = \langle \phi_{k+1}-\phi_k, \phi_{k+1} \rangle - \frac{1}{2}\| \phi_{k+1} - \phi_k \|^2$, we have
\begin{align*}
&\quad \frac{1+\alpha h}{2}\| \phi_{k+1} \|^2 - \frac{1}{2}\| \phi_k \|^2 \\
&= \underbrace{ - \frac{\alpha h}{2}(1+\alpha h)\| x_{k+1} - x_k \|^2 }_{\le 0}
   - \gamma h^2 (1+\alpha h) \langle \nabla f(x_k), x_{k+1}-x_k \rangle \\
&\quad - \left( \frac{1}{2}\gamma^2 + \gamma - \frac{1}{2}\alpha h \right) h^4 \| \nabla f(x_k) \|^2 \\
&\le \underbrace{ - \gamma h^2 (1+\alpha h) \langle \nabla f(x_k), x_{k+1}-x_k \rangle }_{\mathbf{I}}
    - \left( \frac{1}{2}\gamma^2 + \gamma - \frac{1}{2}\alpha h \right) h^4 \| \nabla f(x_k) \|^2.  
\end{align*}
Using Lemma \ref{Lem: func_bound}, we obtain
\begin{align*}
\mathbf{I}
&\le -\gamma h^2 \left[ \psi_k - \psi_{k-1} - \gamma h^2 \| \nabla f(x_k) \|^2 \right],
\end{align*}
so that
$$\frac{1+\alpha h}{2}\| \phi_{k+1} \|^2 - \frac{1}{2}\| \phi_k \|^2 
   + \gamma h^2 (\psi_k - \psi_{k-1}) 
\le - \left[ \gamma(2-\gamma) - \alpha h \right] \frac{h^4}{2} \| \nabla f(x_k) \|^2. $$
Adding $\gamma\alpha h^3 \psi_k = \gamma\alpha h^3\left( f(x_k) - f^* - \frac{h^2}{2}\| \nabla f(x_k) \|^2 \right)$ to both sides, we arrive at
$$ (1+\alpha h)E_{k+1}-E_k 
\le - \left[ \gamma(2-\gamma) + (\gamma-1) \alpha h \right] \frac{h^4}{2} \| \nabla f(x_k) \|^2
    + \gamma\alpha h^3 \left( f(x_k) - f^* \right), $$
as claimed.
\end{proof}

In what follows, we give the main result of this subsection.

\begin{theorem}\label{Thm: algo_PL_OGM}
Let $f:H\to\R$ be convex, $L$-smooth and satisfy the Polyak-\L ojasiewicz inequality \eqref{E: PL} with constant $\mu$ for $L>\mu>0$. Let $(x_k)$ and $(y_k)$ be generated by \eqref{Algo: AGM}, where
$$ h=1/\sqrt{L},\qbox{}\alpha h = \frac{2q}{1+\sqrt{2q-q^2}},\qbox{}\gamma = \frac{\sqrt{2q-q^2}-q}{1-q}, $$
with $q=\mu/L$. If $y_0 = -(\gamma-\alpha h)h^2\nabla f(x_0)+x_0$, then for every $k\ge 1$,
$$ f(y_k) - f^* \le \left( f(x_0) - f^* \right) (1 + \rho)^{-k+1}, $$
where $\rho = \frac{2q}{1+\sqrt{2q-q^2}}$. 
\end{theorem}

\begin{proof}
Using Lemma \ref{Lem: E_k_bound_PL}, together with the Polyak-\L ojasiewicz inequality:
$$ f(x_k) - f^* \le \frac{1}{2\mu}\| \nabla f(x_k) \|^2, $$
we deduce that 
$$ (1+\alpha h)E_{k+1}-E_k 
\le - \underbrace{\left[ \gamma(2-\gamma) + (\gamma-1) \alpha h - \frac{\gamma \alpha}{\mu h} \right] }_{\mathbf{II}} \frac{h^4}{2} \| \nabla f(x_k) \|^2. $$
As one can see, $\mathbf{II}\ge 0$ holds as long as $\gamma\in(0,2)$ and $\alpha h \le \frac{\gamma(2-\gamma)q}{(1-q)\gamma+q}$. It can be verified that
\begin{equation}\nonumber
\gamma = \frac{\sqrt{2q-q^2}-q}{1-q}\in(0,1),\quad
\alpha h = \frac{2q}{1+\sqrt{2q-q^2}} = \frac{\gamma(2-\gamma)q}{(1-q)\gamma+q},
\end{equation}
satisfy the aforementioned condition. As a result, $E_k \le (1+\alpha h)^{-k+1}E_1$. Recalling Remark \ref{Rem: x_0}, it can be shown that the choice of $y_0 = -(\gamma-\alpha h)h^2\nabla f(x_0)+x_0$ gives
$$ \phi_1 = x_1 - y_1 = x_1-x_0+h^2\nabla f(x_0) = 0. $$
Recalling \eqref{E: E_k_PL}, we have $E_1 = \gamma h^2 \psi_0 \le \gamma h^2 \left( f(x_0) - f^* \right)$. It follows that
$$ f(y_k) - f^* \le \psi_{k-1} \le \frac{E_k}{\gamma h^2} \le \frac{E_1}{\gamma h^2}(1+\alpha h)^{-k+1} \le \left( f(x_0) - f^* \right)(1+\alpha h)^{-k+1}, $$
as claimed.
\end{proof}

\begin{remark}
For $q$ small enough, one obtains
$$ f(y_{k+1}) - f^* = \mathcal{O}\left( (1+2q)^{-k} \right). $$
This rate is inferior to the one obtained in Theorem \ref{Thm: algo_QG_OGM}. Although the Polyak-\L ojasiewicz inequality implies a quadratic growth condition \cite{Bolte_2017}, the assumption of a unique minimizer is not required here, which may account for the inferiority of the convergence rate. On the other hand, this rate is comparable to the one obtained by gradient descent for convex functions that satisfies the Polyak-\L ojasiewicz inequality \cite[Corollary 3.1]{Zhang_2019}. An interesting thing is that this rate was derived under a choice $\gamma<1$ instead of $\gamma=1$, whose convergence rate is only $\mathcal{O}\left( (1+q)^{-k} \right)$. This demonstrates one potential advantage of including $\gamma$: one can obtain more control on manipulating the convergence behavior of the algorithms and derive a faster convergence rate. The interested reader can refer to \cite{Wang_2025_MAPM} for an example.    
\end{remark}

\section{Conclusions}\label{Sec: conclusions}

We have provided a dynamic interpretation for the extra speedup enjoyed by state-of-the-art first-order algorithms for smooth strongly convex functions with respect to Nesterov's method. We prove this idea within the framework of inertial dynamics and algorithms, and recover the fastest known rates under strong convexity. Following the same strategy, we improve the convergence rate by a factor of $\sqrt{2}$ for convex functions with quadratic growth and a unique minimizer.

\begin{acknowledgements}
This work was partially funded by the China Scholarship Council~202208520010, and also benefited from the support of the FMJH Program Gaspard Monge for optimization and operations research and their interactions with data science.
\end{acknowledgements}

\bibliography{myrefs}

\begin{appendices}

\section{Geometric properties} \label{A:geometry}

In this section, we recall some terminology concerning the geometric assumptions that will be made on $f$ to establish different convergence rates. For a general account on convexity, the reader is referred to \cite{Bauschke_2017,Peypouquet_2015}.

In all that follows $f:H\to\R$ is differentiable and attains its minimum $f^*$ at $S=\argmin(f)$.

The function $f$ is {\it $L$-smooth} if its gradient $\nabla f$ is Lipschitz-continuous with constant $L$. Smooth convex functions can be characterized as follows: $f$ is $L$-smooth and convex if, and only if, 
$$f(y)\ge f(x)+\langle \nabla f(x),y-x\rangle+\frac{1}{2L}\|\nabla f(x)-\nabla f(y)\|^2$$
for every $x,y\in H$.

Now, $f$ is {\it $\mu$-strongly convex} if
$$f(y)\ge f(x)+\langle \nabla f(x),y-x\rangle+\frac{\mu}{2}\|x-y\|^2$$
for every $x,y\in H$. If $f$ is strongly convex, it has a unique minimizer $x^*$ and 
$$f(x) - f^*\ge \frac{\mu}{2}\|x-x^*\|^2.$$
More generally, $f$ has {\it $\mu$-quadratic growth} if 
$$f(x) - f^*\ge \frac{\mu}{2}\dist(x,S)^2.$$
On the other hand, $f$ satisfies the {\it Polyak-\L ojasiewicz} inequality with constant $\mu > 0$, or is a {\it $\mu$-P\L\ function}, if
$$\| \nabla f(x) \|^2 \ge 2\mu \big( f(x) - f^* \big),$$
for every $x\in H$. The interested reader is referred to \cite{Bolte_2017,Polyak_1963,Lojasiewicz_1963,Lojasiewicz_1958,Lojasiewicz_1959}. If $f$ is $\mu$-strongly convex, by Young's inequality, we have
$$ \frac{\mu}{2}\| x - x^* \|^2 + \frac{1}{2\mu}\| \nabla f(x) \|^2 \ge \langle \nabla f(x),x-x^*\rangle \ge f(x)-f^*+\frac{\mu}{2}\|x-x^*\|^2,$$
and so $f$ is $\mu$-P\L. It follows that both quadratic growth and the P{\L} inequality are generalized forms of strong convexity. For the readers' convenience, we include a short discussion on the relationship between these two notions, mostly adapted from 
\cite{Bolte_2017} (see also \cite{Bolte_2007,Bolte_2010}).

\begin{lemma}
Let $f:H\to\R$ be smooth, and such that
\begin{equation} \label{E:desingularizing}
\varphi'\big(f(x)-f^*\big)\|\nabla f(x)\|\ge 1    
\end{equation}
for all $x\in H$, where $\varphi$ is strictly increasing and differentiable, with $\varphi(0)=0$. Let $u:[0,\infty)\to H$ satisfy $-\dot u(t)=\nabla f\big(u(t)\big)$ for $t>0$. For $0\le s\le t$, we have  
\begin{equation} \label{E:length}
\|u(t)-u(s)\|\le \varphi\big(f(u(s))-f^*\big)-\varphi\big(f(u(t))-f^*\big),  
\end{equation}
and
\begin{equation} \label{E:Loja_values}
\frac{d}{dt}\left[\frac{}{}\!f(u(t))-f^*\right]\le -\left[\frac{1}{\varphi'\big(f(u(t))-f^*\big)}\right]^2.   
\end{equation}
\end{lemma}

\begin{proof}
Let $0\le s\le t$. We have
\begin{eqnarray*}
\varphi\big(f(u(s))-f^*\big)-\varphi\big(f(u(t))-f^*\big) 
& = & \int_t^s\frac{d}{d\tau} \big(\varphi\circ(f-f^*)\big)\big(u(\tau)\big) \,d\tau \\
& = & \int_t^s\varphi'\big(f(u(\tau)-f^*)\big)\left\langle \nabla f\big(u(\tau)\big),\dot u(\tau)\right\rangle\,d\tau  \\
& = & \int_s^t\varphi'\big(f(u(\tau)-f^*\big)\left\|\nabla f\big(u(\tau)\big)\right\|\left\|\dot u(\tau)\right\|\,d\tau  \\
& \ge & \int_s^t\left\|\dot u(\tau)\right\|\,d\tau .
\end{eqnarray*}
Therefore,
$$\|u(t)-u(s)\|= \left\|\int_s^t\dot u(\tau)\,d\tau\right\|\le \int_s^t\left\|\dot u(\tau)\right\|\,d\tau\le \varphi\big(f(u(s))-f^*\big)-\varphi\big(f(u(t))-f^*\big).$$  
On the other hand, we have
$$\frac{d}{dt}\left[\frac{}{}\!f(u(t))-f^*\right]=\langle \nabla f\big(u(t)\big),\dot u(t)\rangle=-\|\nabla f\big(u(t)\big)\|^2\le -\left[\frac{1}{\varphi'\big(f(u(t))-f^*\big)}\right
]^2,$$
as required.
\end{proof}

If $f$ is $\mu$-P{\L}, it satisfies \eqref{E:desingularizing} with $\varphi=\sqrt{\frac{2s}{\mu}}$.

\begin{corollary} \label{C:PLtoQG}
If $f$ is smooth and $\mu$-P{\L}, it has $\mu$-quadratic growth.
\end{corollary}

\begin{proof}
Let $\varphi(s)=\sqrt{\frac{2s}{\mu}}$, so that $\varphi'(s)=\frac{1}{\sqrt{2\mu s}}$. Let $u:[0,\infty)\to H$ be the unique solution of $-\dot u(t)=\nabla f\big(u(t)\big)$, with $u(0)=z$. From \eqref{E:Loja_values}, we obtain
$$\frac{d}{dt}\left[\frac{}{}\!f(u(t))-f^*\right]\le -2\mu \big(f(u(t))-f^*\big),$$
and so 
$$f(u(t))-f^*\le e^{-2\mu t}\big(f(z)-f^*\big),$$
which tends to zero as $t\to\infty$. By \eqref{E:length}, $u(t)$ converges strongly, as $t\to\infty$. Its limit, say $\hat u$, must be a minimizer of $f$. In particular,
$$\dist(z,S)\le\|\hat u-z\|\le \varphi\big( f(z)-f^*\big)=\sqrt{\frac{2(f(z)-f^*)}{\mu}},$$
as claimed.
\end{proof}

Conversely, we have the following:
\begin{proposition} \label{P:QGtoPL}
    If $f$ is convex and has $\mu$-quadratic growth, it is $\frac{\mu}{4}$-P{\L}.
\end{proposition}

\begin{proof}
Let $f$ be convex, and let $x_S$ be the projection of $x$ onto $S$. We have
$$f(x)-f^*\le\langle \nabla f(x),x-x_S\rangle\le \|\nabla f(x)\|\dist(x,S).$$
By $\mu$-quadratic growth, we have
$$f(x)-f^*\le \sqrt{\frac{2}{\mu}}\|\nabla f(x)\|\sqrt{f(x)-f^*},$$
which we can rewrite as
$$\| \nabla f(x) \|^2 \ge 2\left[\frac{\mu}{4}\right] \big( f(x) - f^* \big)$$
to conclude.
\end{proof}

\section{Proof of Lemma \ref{Lem: W_bound}} \label{A:Lemma_continuous}

\begin{proof}
We begin by using \eqref{System: HBF-H}, to compute the time derivative of $W$:
\begin{equation}\label{E: W_dot}
\begin{aligned}
\dot{W} 
&= \big\langle \dot{x} + \beta\nabla f(x) + \xi(x-x^*), \ddot{x} + \beta\nabla^2 f(x) \dot{x} + \xi\dot{x} \big\rangle \\ 
&\quad - \omega\xi(\alpha-\xi)\langle x-x^*, \dot{x} \rangle 
  + \theta \langle \nabla f(x), \dot{x} \rangle \\
&= \big\langle \dot{x} + \beta\nabla f(x) + \xi(x-x^*), -(\alpha-\xi)\dot{x} - \gamma\nabla f(x) \big\rangle \\
&\quad - \omega\xi(\alpha-\xi)\langle x-x^*, \dot{x} \rangle
  + \theta \langle \nabla f(x), \dot{x} \rangle \\
&= -(\alpha-\xi)\|\dot{x}\|^2 
   - [\gamma + (\alpha-\xi)\beta - \theta]\langle \nabla f(x),\dot{x} \rangle
   - \beta\gamma\|\nabla f(x)\|^2 \\
&\quad \underbrace{ - (1+\omega)\xi(\alpha-\xi) \langle x-x^*,\dot{x} \rangle }
   - \xi\gamma\langle\nabla f(x),x-x^* \rangle.
\end{aligned}
\end{equation}
Then, by developing \eqref{E: W}, we obtain
\begin{align*}
W 
&= \frac{1}{2}\|\dot{x}\|^2 
  + \frac{1}{2}\beta^2\|\nabla f(x) \|^2 
  + \frac{1}{2}\xi[\xi-\omega(\alpha-\xi)]\|x-x^*\|^2
  + \theta(f(x)-f^*) \\
&\quad + \beta\langle \nabla f(x),\dot{x} \rangle
  + \underbrace{ \xi\langle x-x^*,\dot{x} \rangle }
  + \xi\beta\langle \nabla f(x),x-x^*\rangle.
\end{align*}
To cancel the term $\langle x-x^*,\dot{x} \rangle$ in \eqref{E: W_dot}, we multiply this equation by $(1+\omega)(\alpha-\xi)$, to get
\begin{align*}
(1+\omega)(\alpha-\xi)W 
&= \frac{1}{2}(1+\omega)(\alpha-\xi)\|\dot{x}\|^2 
  + \frac{1}{2}\beta^2(1+\omega)(\alpha-\xi) \|\nabla f(x) \|^2 \\
&\quad + \frac{1}{2}(1+\omega)(\alpha-\xi)\xi[\xi-\omega(\alpha-\xi)]\|x-x^*\|^2 \\
&\quad + \theta(1+\omega)(\alpha-\xi)(f(x)-f^*) 
  + (1+\omega)(\alpha-\xi)\beta\langle \nabla f(x),\dot{x} \rangle \\
&\quad + \underbrace{ (1+\omega)\xi(\alpha-\xi)\langle x-x^*,\dot{x} \rangle }
  + (1+\omega)\xi(\alpha-\xi)\beta\langle \nabla f(x),x-x^*\rangle,
\end{align*}
and add it to \eqref{E: W_dot} so that
\begin{align*}
 \dot{W} + (1+\omega)(\alpha-\xi)W 
&= -\frac{1}{2}(1-\omega)(\alpha-\xi)\|\dot{x}\|^2
   - [\gamma - \omega(\alpha-\xi)\beta - \theta]\langle \nabla f(x),\dot{x} \rangle \\
&\quad - \beta\left[ \gamma - \frac{1}{2}(1+\omega)(\alpha-\xi)\beta \right]\|\nabla f(x)\|^2 \\
&\quad + \frac{1}{2}(1+\omega)(\alpha-\xi)\xi[\xi-\omega(\alpha-\xi)]\|x-x^*\|^2 \\
&\quad + \theta(1+\omega)(\alpha-\xi)(f(x)-f^*) \\
&\quad - \xi[\gamma-(1+\omega)(\alpha-\xi)\beta]\langle\nabla f(x),x-x^* \rangle.
\end{align*}
Substituting $\gamma = \theta + (\alpha-\xi)\beta$, we arrive at
\begin{align*}
 \dot{W} + (1+\omega)(\alpha-\xi)W 
&= \underbrace{ -\frac{1}{2}(1-\omega)(\alpha-\xi)\|\dot{x}\|^2
   - (1-\omega)(\alpha-\xi)\beta \langle \nabla f(x),\dot{x} \rangle }_{\mathbf{I}} \\
&\quad - \beta\left[ \theta + \frac{1}{2}(1-\omega)(\alpha-\xi)\beta \right]\|\nabla f(x)\|^2 \\
&\quad + \frac{1}{2}(1+\omega)(\alpha-\xi)\xi[\xi-\omega(\alpha-\xi)]\|x-x^*\|^2 \\
&\quad + \theta(1+\omega)(\alpha-\xi)(f(x)-f^*)
  \underbrace{ - \xi[\theta-\omega(\alpha-\xi)\beta]\langle\nabla f(x),x-x^* \rangle }_{\mathbf{II}}.
\end{align*}
Using
\begin{align*}
\mathbf{I}
&= -\frac{1}{2}(1-\omega)(\alpha-\xi) \left\| \dot{x} + \beta\nabla f(x) \right\|^2
   + \frac{1}{2}(1-\omega)(\alpha-\xi)\beta^2 \| \nabla f(x) \|^2 \\
&\le \frac{1}{2}(1-\omega)(\alpha-\xi)\beta^2 \| \nabla f(x) \|^2,  
\end{align*}
and 
\begin{align*}
\mathbf{II}
&= - \xi\theta \langle\nabla f(x),x-x^* \rangle
   + \omega\xi(\alpha-\xi)\beta \langle\nabla f(x),x-x^* \rangle \\
&= - \xi\theta \langle\nabla f(x),x-x^* \rangle
   + \frac{1}{2}\omega\beta \| \nabla f(x) \|^2
   + \frac{1}{2}\omega\beta \xi^2(\alpha-\xi)^2 \| x-x^* \|^2 \\
&\quad \underbrace{ -\frac{1}{2}\omega\beta \left\| \nabla f(x) - \xi(\alpha-\xi)(x-x^*) \right\|^2 }_{\le 0} \\
&\le - \xi\theta \langle\nabla f(x),x-x^* \rangle
   + \frac{1}{2}\omega\beta \| \nabla f(x) \|^2
   + \frac{1}{2}\omega\beta \xi^2(\alpha-\xi)^2 \| x-x^* \|^2, 
\end{align*}
we deduce that
\begin{align*}
&\quad \dot{W} + (1+\omega)(\alpha-\xi)W \\ 
&\le - \beta\left( \theta - \frac{\omega}{2} \right)\|\nabla f(x)\|^2 
  + (1+\omega)(\alpha-\xi)\theta (f(x)-f^*) 
  - \xi\theta \langle\nabla f(x),x-x^* \rangle\\
&\quad + \frac{1}{2}(\alpha-\xi)\xi[ (1+\omega)(\xi-\omega(\alpha-\xi)) + \omega\xi(\alpha-\xi)\beta ]\|x-x^*\|^2,
\end{align*}
which allows us to conclude. 
\end{proof}

\section{Proof of Lemma \ref{Lem: E_k_bound_SC}} \label{Proof: Lem: E_k_bound_SC}

\begin{proof}
By definition of $\phi_k$, we have
\begin{align*}
\phi_{k+1} 
&= (x_{k+1}-y_{k+1}) + \frac{ \frac{1}{2}\alpha h }{ 1 + \frac{1}{2}\alpha h }(x_{k+1}-x^*) \\
&= (x_{k+1}-x_k) + \frac{ \frac{1}{2}\alpha h }{ 1 + \frac{1}{2}\alpha h }(x_{k+1}-x^*) + h^2 \nabla f(x_k)\\
&= \frac{ 1 + \alpha h }{ 1 + \frac{1}{2}\alpha h }(x_{k+1}-x_k)
  + \frac{ \frac{1}{2}\alpha h }{ 1 + \frac{1}{2}\alpha h }(x_k-x^*) + h^2 \nabla f(x_k), 
\end{align*}
and using \eqref{E: iterates},
\begin{align*}
\phi_{k+1}-\phi_k
&= -\frac{\alpha h}{2}\left( \frac{ 1 + \alpha h }{ 1 + \frac{1}{2}\alpha h } \right) (x_{k+1}-x_k) - \gamma h^2 \nabla f(x_k) \\
&= -\frac{\alpha h}{4}\phi_{k+1}
   - \frac{\alpha h}{4}\left( \frac{ 1 + \alpha h }{ 1 + \frac{1}{2}\alpha h } \right)(x_{k+1}-x_k)
   + \frac{\frac{1}{8}\alpha^2 h^2}{1+\frac{1}{2}\alpha h}(x_k-x^*) \\
&\quad -\left( \gamma - \frac{1}{4}\alpha h \right) h^2 \nabla f(x_k).
\end{align*}
It follows that
\begin{align*}
\| \phi_{k+1}-\phi_k \|^2
&= \frac{\alpha^2 h^2}{4}\left( \frac{ 1 + \alpha h }{ 1 + \frac{1}{2}\alpha h } \right)^2\| x_{k+1}-x_k \|^2 
  + \gamma^2 h^4 \| \nabla f(x_k) \|^2 \\
&\quad + \gamma\alpha h^3 \left( \frac{ 1 + \alpha h }{ 1 + \frac{1}{2}\alpha h } \right) \langle \nabla f(x_k), x_{k+1} - x_k \rangle,
\end{align*}
and
\begin{align*}
&\quad \langle \phi_{k+1}-\phi_k, \phi_{k+1} \rangle \\
&= - \frac{\alpha h}{4}\| \phi_{k+1} \|^2
   - \frac{\alpha h}{4} \left(\frac{ 1 + \alpha h }{ 1 + \frac{1}{2}\alpha h }\right)^2\| x_{k+1}-x_k \|^2 \\
&\quad - \left( \gamma - \frac{1}{2}\alpha h \right) h^2 \left(\frac{ 1 + \alpha h }{ 1 + \frac{1}{2}\alpha h }\right) \langle \nabla f(x_k),x_{k+1}-x_k \rangle
   + \frac{\frac{1}{16}\alpha^3 h^3 \| x_k - x^* \|^2 }{\left( 1+\frac{1}{2}\alpha h \right)^2} \\
&\quad - \left( \gamma - \frac{1}{2}\alpha h \right)\frac{ \frac{1}{2}\alpha h^3 }{1+\frac{1}{2}\alpha h}\langle \nabla f(x_k), x_k-x^* \rangle
   -\left( \gamma - \frac{1}{4}\alpha h \right) h^4 \| \nabla f(x_k) \|^2.
\end{align*}
Since $\frac{1}{2}\| \phi_{k+1} \|^2 - \frac{1}{2}\| \phi_k \|^2 = \langle \phi_{k+1}-\phi_k, \phi_{k+1} \rangle - \frac{1}{2}\| \phi_{k+1} - \phi_k \|^2$, we have
\begin{align*}
&\quad \left( 1 + \frac{1}{2}\alpha h \right)\frac{1}{2}\| \phi_{k+1} \|^2 - \frac{1}{2}\| \phi_k \|^2\\
&= \underbrace{ - \frac{\alpha h}{4} \left[\frac{ (1 + \alpha h)^2 }{ 1 + \frac{1}{2}\alpha h }\right]\| x_{k+1}-x_k \|^2 }_{\mathbf{I}} 
  -\left( \frac{\gamma^2}{2} + \gamma - \frac{1}{4}\alpha h \right) h^4 \| \nabla f(x_k) \|^2\\
&\quad - \left[ \gamma \left( 1 + \frac{1}{2}\alpha h \right) - \frac{1}{2}\alpha h \right] h^2 \left(\frac{ 1 + \alpha h }{ 1 + \frac{1}{2}\alpha h }\right) \langle \nabla f(x_k),x_{k+1}-x_k \rangle \\
&\quad + \frac{\frac{1}{16}\alpha^3 h^3}{\left( 1+\frac{1}{2}\alpha h \right)^2}\| x_k - x^* \|^2 
 - \left( \gamma - \frac{1}{2}\alpha h \right)\frac{ \frac{1}{2}\alpha h^3 }{1+\frac{1}{2}\alpha h}\langle \nabla f(x_k), x_k-x^* \rangle,
\end{align*}
which, together with
\begin{align*}
\mathbf{I}
&= - \frac{\alpha h}{4} \left[\frac{ (1 + \alpha h)^2 }{ 1 + \frac{1}{2}\alpha h }\right] \left\| x_{k+1}-x_k + \frac{\gamma h^2}{1+\alpha h} \nabla f(x_k) \right\|^2 \\
&\quad + \frac{\gamma\alpha h^3}{2} \left( \frac{1+\alpha h}{1+\frac{1}{2}\alpha h} \right) \langle \nabla f(x_k), x_{k+1}-x_k \rangle
   + \frac{ \frac{1}{4}\gamma^2 \alpha h }{1 + \frac{1}{2}\alpha h}h^4 \| \nabla f(x_k) \|^2\\
&\le \frac{\gamma\alpha h^3}{2} \left( \frac{1+\alpha h}{1+\frac{1}{2}\alpha h} \right) \langle \nabla f(x_k), x_{k+1}-x_k \rangle
   + \frac{ \frac{1}{4}\gamma^2 \alpha h }{1 + \frac{1}{2}\alpha h}h^4 \| \nabla f(x_k) \|^2,    
\end{align*}
gives
\begin{align*}
&\quad \left( 1 + \frac{1}{2}\alpha h \right)\frac{1}{2}\| \phi_{k+1} \|^2 - \frac{1}{2}\| \phi_k \|^2\\
&\le \underbrace{ - \left( \gamma - \frac{1}{2}\alpha h \right) h^2 \left(\frac{ 1 + \alpha h }{ 1 + \frac{1}{2}\alpha h }\right) \langle \nabla f(x_k),x_{k+1}-x_k \rangle }_{\mathbf{II}}
   + \frac{\frac{1}{16}\alpha^3 h^3 \| x_k - x^* \|^2 }{\left( 1+\frac{1}{2}\alpha h \right)^2} \\
&\quad \underbrace{ - \left( \gamma - \frac{1}{2}\alpha h \right)\frac{ \frac{1}{2}\alpha h^3 }{1+\frac{1}{2}\alpha h}\langle \nabla f(x_k), x_k-x^* \rangle }_{\mathbf{III}} \\
&\quad -\left[ \frac{\gamma^2}{2} + \gamma - \frac{1}{4}\alpha h \left( 1 + \frac{\gamma^2}{1+\frac{1}{2}\alpha h} \right) \right] h^4 \| \nabla f(x_k) \|^2.
\end{align*}
Recalling that
$$ \theta = \frac{ \gamma-\frac{1}{2}\alpha h }{1 + \frac{1}{2}\alpha h}, $$
we use Lemma \ref{Lem: func_bound}, to obtain
\begin{align*}
\mathbf{II}
&\le - \theta h^2 (\psi_k-\psi_{k-1}) + \frac{ \gamma\left(\gamma-\frac{1}{2}\alpha h \right)}{1 + \frac{1}{2}\alpha h} h^4 \| \nabla f(x_k) \|^2 \\
&= - \theta h^2 (\psi_k-\psi_{k-1}) + \gamma\left(\gamma-\frac{1}{2}\alpha h \right)\left( 1 - \frac{ \frac{1}{2}\alpha h }{1 + \frac{1}{2}\alpha h} \right) h^4 \| \nabla f(x_k) \|^2, 
\end{align*}
and
\begin{align*}
\mathbf{III}
&\le - \left( \frac{1}{2}\alpha h \right)\theta h^2 \psi_k 
  - \left( \gamma - \frac{1}{2}\alpha h \right)\frac{ \frac{1}{4}\alpha h^5 }{1+\frac{1}{2}\alpha h}\| \nabla f(x_k) \|^2 \\ 
&\quad - \left( \gamma - \frac{1}{2}\alpha h \right)\frac{ \frac{1}{4}\mu\alpha h^3 }{1+\frac{1}{2}\alpha h}\| x_k-x^* \|^2 \\
&\le 
  - \left( \frac{1}{2}\alpha h \right)\theta h^2 \psi_k 
  - \left( \gamma - \frac{1}{2}\alpha h \right)\frac{ \frac{1}{4}\alpha h^5 }{1+\frac{1}{2}\alpha h}\| \nabla f(x_k) \|^2 \\ 
&\quad - \frac{ \frac{1}{4}\mu\gamma\alpha h^3 }{\left(1+\frac{1}{2}\alpha h \right)^2}\| x_k-x^* \|^2 
   \quad( \text{since }\gamma\ge1+\alpha h ).
\end{align*}
Whence, we deduce that
$$ \left( 1 + \frac{1}{2}\alpha h \right)E_{k+1} - E_k
\le \frac{\frac{1}{16}\alpha h^3 (\alpha^2-4\mu\gamma)}{\left( 1+\frac{1}{2}\alpha h \right)^2}\| x_k - x^* \|^2 - Rh^4\| \nabla f(x_k) \|^2,  $$
where
$$ R= \underbrace{ \frac{\gamma}{2}(2+\alpha h-\gamma) }_{\ge 0} + \frac{\frac{1}{4}\alpha h}{1+\frac{1}{2}\alpha h}\big[ \gamma \underbrace{ (\gamma-\alpha h) }_{\ge 1} + \underbrace{ \gamma - (1+\alpha h) }_{\ge 0} \big] 
\ge \frac{\frac{1}{4}\gamma \alpha h}{1+\frac{1}{2}\alpha h}
\ge \frac{1}{4}\alpha h, $$
since $\gamma\in[1+\alpha h, 2+\alpha h]$. The result follows by using $R\ge \frac{1}{4}\alpha h$.
\end{proof}

\section{Proof of Lemma \ref{Lem: E_k_bound_SC_TMM}}
\label{Proof: Lem: E_k_bound_SC_TMM}

\begin{proof}
By definition of $\phi_k$, we have
\begin{align*}
\phi_{k+1} 
&= (x_{k+1}-y_{k+1}) + \frac{2}{3}\alpha h (x_{k+1}-x^*) \\
&= (x_{k+1}-x_k) + \frac{2}{3}\alpha h (x_{k+1}-x^*) + h^2 \nabla f(x_k)\\
&= \left( 1 + \frac{2}{3}\alpha h \right)(x_{k+1}-x_k)
  + \frac{2}{3}\alpha h (x_k-x^*) + h^2 \nabla f(x_k), 
\end{align*}
and using \eqref{E: iterates},
\begin{align*}
\phi_{k+1}-\phi_k
&= -\frac{\alpha h}{3} (x_{k+1}-x_k) - \gamma h^2 \nabla f(x_k) \\
&= -\frac{ \frac{1}{3}\alpha h }{ 1 + \frac{2}{3}\alpha h }\phi_{k+1}
   + \frac{\frac{2}{9}\alpha^2 h^2}{1+\frac{2}{3}\alpha h}(x_k-x^*)  
   - \left( \gamma - \frac{ \frac{1}{3}\alpha h }{ 1 + \frac{2}{3}\alpha h } \right) h^2 \nabla f(x_k).
\end{align*}
As a result,
$$ \| \phi_{k+1}-\phi_k \|^2
= \frac{\alpha^2 h^2}{9}\| x_{k+1}-x_k \|^2 
  + \gamma^2 h^4 \| \nabla f(x_k) \|^2 
  + \frac{2\gamma\alpha h^3}{3} \langle \nabla f(x_k), x_{k+1} - x_k \rangle, $$
and
\begin{align*}
&\quad \langle \phi_{k+1}-\phi_k, \phi_{k+1} \rangle \\
&= - \frac{ \frac{1}{3}\alpha h }{ 1 + \frac{2}{3}\alpha h }\| \phi_{k+1} \|^2 
   + \frac{2}{9}\alpha^2 h^2\langle x_{k+1}-x_k, x_k - x^* \rangle\\
&\quad + \frac{\frac{4}{27}\alpha^3 h^3}{1+\frac{2}{3}\alpha h}\| x_k - x^* \|^2 
   - \left( \gamma - \frac{\frac{2}{3}\alpha h}{1+\frac{2}{3}\alpha h} \right)\frac{2}{3}\alpha h^3\langle \nabla f(x_k), x_k-x^* \rangle \\
&\quad - \left[ \gamma \left( 1 + \frac{2}{3}\alpha h \right) - \frac{1}{3}\alpha h \right] h^2  \langle \nabla f(x_k),x_{k+1}-x_k \rangle 
 - \left( \gamma - \frac{ \frac{1}{3}\alpha h }{ 1 + \frac{2}{3}\alpha h } \right) h^4 \| \nabla f(x_k) \|^2.
\end{align*}
Since $\frac{1}{2}\| \phi_{k+1} \|^2 - \frac{1}{2}\| \phi_k \|^2 = \langle \phi_{k+1}-\phi_k, \phi_{k+1} \rangle - \frac{1}{2}\| \phi_{k+1} - \phi_k \|^2$, we have
\begin{equation}\label{E: phi_k_diff_SC_TMM}
\begin{aligned}
&\quad \left( 1 + \frac{\frac{2}{3}\alpha h}{1 + \frac{2}{3}\alpha h} \right)\frac{1}{2}\| \phi_{k+1} \|^2 - \frac{1}{2}\| \phi_k \|^2\\
&= \underbrace{ \frac{2}{9}\alpha^2 h^2\langle x_{k+1}-x_k, x_k - x^* \rangle }
   + \frac{\frac{4}{27}\alpha^3 h^3}{1+\frac{2}{3}\alpha h}\| x_k - x^* \|^2 \\
&\quad - \frac{\alpha^2 h^2}{18}\| x_{k+1}-x_k \|^2 
   - \left( \gamma - \frac{\frac{2}{3}\alpha h}{1+\frac{2}{3}\alpha h} \right)\frac{2}{3}\alpha h^3\langle \nabla f(x_k), x_k-x^* \rangle \\
&\quad - \left[ \gamma ( 1 + \alpha h ) - \frac{1}{3}\alpha h \right] h^2  \langle \nabla f(x_k),x_{k+1}-x_k \rangle 
 - \left( \frac{\gamma^2}{2} + \gamma - \frac{ \frac{1}{3}\alpha h }{ 1 + \frac{2}{3}\alpha h } \right) h^4 \| \nabla f(x_k) \|^2.
\end{aligned}
\end{equation}
Likewise, using \eqref{E: iterates} and recalling that
$$ y_{k+1}-x^* = x_k - x^* - h^2\nabla f(x_k), $$
we have
\begin{align*}
y_{k+1} - y_k
&= (1+\alpha h)(x_{k+1}-x_k) + \gamma h^2 \nabla f(x_k) \\
&= -\frac{\frac{1}{3}\alpha h}{1 + \frac{2}{3}\alpha h}(y_{k+1}-x^*)
   + (1+\alpha h)(x_{k+1}-x_k) \\ 
&\quad + \frac{\frac{1}{3}\alpha h}{1 + \frac{2}{3}\alpha h}(x_k - x^*) 
   + \left( \gamma - \frac{\frac{1}{3}\alpha h}{1 + \frac{2}{3}\alpha h} \right) h^2 \nabla f(x_k), 
\end{align*}
so that
\begin{align*}
&\quad \langle y_{k+1}-y_k, y_{k+1}-x^* \rangle \\
&= -\frac{\frac{1}{3}\alpha h}{1 + \frac{2}{3}\alpha h}\| y_{k+1}-x^* \|^2
  + (1+\alpha h)\langle x_{k+1}-x_k, x_k - x^* \rangle\\
&\quad - (1+\alpha h)h^2 \langle \nabla f(x_k), x_{k+1}-x_k \rangle
  + \frac{\frac{1}{3}\alpha h}{1 + \frac{2}{3}\alpha h}\| x_k - x^* \|^2 \\
&\quad + \left( \gamma - \frac{\frac{2}{3}\alpha h}{1 + \frac{2}{3}\alpha h} \right) h^2 \langle \nabla f(x_k), x_k - x^* \rangle 
  - \left( \gamma - \frac{\frac{1}{3}\alpha h}{1 + \frac{2}{3}\alpha h} \right) h^4 \| \nabla f(x_k) \|^2.
\end{align*}
Since $\frac{1}{2}\| y_{k+1}-x^* \|^2 - \frac{1}{2}\| y_k-x^* \|^2 = \langle y_{k+1}-y_k, y_{k+1}-x^* \rangle - \frac{1}{2}\| y_{k+1} - y_k \|^2$, we have
\begin{align*}
&\quad \left( 1 + \frac{\frac{2}{3}\alpha h}{1 + \frac{2}{3}\alpha h} \right)\frac{1}{2}\| y_{k+1}-x^* \|^2 - \frac{1}{2}\| y_k-x^* \|^2 \\
&= \underbrace{ (1+\alpha h)\langle x_{k+1}-x_k, x_k - x^* \rangle }
  - (1+\alpha h)h^2 \langle \nabla f(x_k), x_{k+1}-x_k \rangle \\
&\quad + \frac{\frac{1}{3}\alpha h}{1 + \frac{2}{3}\alpha h}\| x_k - x^* \|^2 
  + \left( \gamma - \frac{\frac{2}{3}\alpha h}{1 + \frac{2}{3}\alpha h} \right) h^2 \langle \nabla f(x_k), x_k - x^* \rangle \\ 
&\quad - \left( \gamma - \frac{\frac{1}{3}\alpha h}{1 + \frac{2}{3}\alpha h} \right) h^4 \| \nabla f(x_k) \|^2 
  - \frac{1}{2}\| y_{k+1} - y_k \|^2.
\end{align*}
To cancel the term $\langle x_{k+1}-x_k, x_k - x^* \rangle$ in \eqref{E: phi_k_diff_SC_TMM}, we multiply this inequality by 
$$ \eta = \frac{\frac{2}{9}\alpha^2 h^2}{1 + \alpha h } $$
and subtract it from \eqref{E: phi_k_diff_SC_TMM}, to get
\begin{align*}
&\quad \left( 1 + \frac{\frac{2}{3}\alpha h}{1 + \frac{2}{3}\alpha h} \right)
  \left( \frac{1}{2}\| \phi_{k+1} \|^2 - \frac{\eta}{2}\| y_{k+1}-x^* \|^2 \right) 
  - \left( \frac{1}{2}\| \phi_k \|^2 - \frac{\eta}{2}\| y_k-x^* \|^2 \right)\\
&= \frac{\frac{2}{27}\alpha^3 h^3 (1+2\alpha h)}{\left( 1+\frac{2}{3}\alpha h \right)(1+\alpha h)}\| x_k - x^* \|^2
  - \frac{\alpha^2 h^2}{18}\| x_{k+1}-x_k \|^2
  + \frac{\frac{1}{9}\alpha^2 h^2}{1 + \alpha h } \| y_{k+1} - y_k \|^2 \\
&\quad \underbrace{ - \left[ \left( 1+\frac{2}{3}\alpha h \right) \gamma - \frac{2}{3}\alpha h \right]\left( \frac{ 1 + \frac{4}{3}\alpha h }{1+\alpha h} \right) \frac{ \frac{2}{3}\alpha h^3 }{1+\frac{2}{3}\alpha h} \langle \nabla f(x_k), x_k-x^* \rangle }_{\mathbf{I}} \\
&\quad \underbrace{ - \left[ \gamma  - \frac{1}{3}\alpha h\left( \frac{1 + \frac{2}{3}\alpha h}{1+\alpha h} \right) \right] ( 1 + \alpha h ) h^2  \langle \nabla f(x_k),x_{k+1}-x_k \rangle }_{\mathbf{II}} \\
&\quad - \left[ \frac{\gamma^2}{2} + \gamma - \frac{ \frac{1}{3}\alpha h }{ 1 + \frac{2}{3}\alpha h } - \left( \gamma - \frac{\frac{1}{3}\alpha h}{1 + \frac{2}{3}\alpha h} \right)\frac{\frac{2}{9}\alpha^2 h^2}{1 + \alpha h } \right] h^4 \| \nabla f(x_k) \|^2.
\end{align*}
By strong convexity, we have
\begin{align*}
\mathbf{I}
&\le - \left[ \left( 1+\frac{2}{3}\alpha h \right) \gamma - \frac{2}{3}\alpha h \right]\left( \frac{ 1 + \frac{4}{3}\alpha h }{1+\alpha h} \right) \frac{ \frac{2}{3}\alpha h^3 }{1+\frac{2}{3}\alpha h} (f(x_k)-f^*) \\
&\quad - \left[ \left( 1+\frac{2}{3}\alpha h \right) \gamma - \frac{2}{3}\alpha h \right]\left( \frac{ 1 + \frac{4}{3}\alpha h }{1+\alpha h} \right) \frac{ \frac{1}{3}\mu\alpha h^3 }{1+\frac{2}{3}\alpha h}\| x_k - x^* \|^2 \\
&\stackrel{\gamma\ge 1}{\le} 
  - \left[ \gamma  - \frac{1}{3}\alpha h\left( \frac{1 + \frac{2}{3}\alpha h}{1+\alpha h} \right) \right] \frac{ \frac{2}{3}\alpha h^3 }{1+\frac{2}{3}\alpha h}(f(x_k)-f^*) \\
&\quad - \frac{ \frac{4}{9}\gamma\alpha^2 h^4 }{(1+\alpha h)\left(1+\frac{2}{3}\alpha h\right)}(f(x_k)-f^*)
  - \frac{ \frac{1}{3}\mu\gamma \alpha h^3 \left( 1 + \frac{4}{3}\alpha h \right) }{(1+\alpha h)\left(1+\frac{2}{3}\alpha h\right)} \| x_k - x^* \|^2.   
\end{align*}
Recalling that
$$ \theta = \gamma  - \frac{1}{3}\alpha h\left( \frac{1 + \frac{2}{3}\alpha h}{1+\alpha h} \right), $$
and using strong convexity:
$$ - \frac{ \frac{4}{9}\gamma\alpha^2 h^4 }{(1+\alpha h)\left(1+\frac{2}{3}\alpha h\right)}(f(x_k)-f^*) \le - \frac{ \frac{2}{9}\mu\gamma\alpha^2 h^4 }{(1+\alpha h)\left(1+\frac{2}{3}\alpha h\right)}\| x_k - x^* \|^2,  $$
we arrive at
\begin{align*}
\mathbf{I}
&\le -\frac{ \frac{2}{3}\alpha h }{1+\frac{2}{3}\alpha h}\theta h^2 \psi_k
   - \left[ \gamma  - \frac{1}{3}\alpha h\left( \frac{1 + \frac{2}{3}\alpha h}{1+\alpha h} \right) \right] \frac{ \frac{1}{3}\alpha h^5 }{1+\frac{2}{3}\alpha h}\| \nabla f(x_k) \|^2 \\
&\quad - \frac{ \frac{1}{3}\mu\gamma \alpha h^3 ( 1 + 2\alpha h ) }{(1+\alpha h)\left(1+\frac{2}{3}\alpha h\right)} \| x_k - x^* \|^2.  
\end{align*}
Likewise, we use Lemma \ref{Lem: func_bound}, to obtain
\begin{align*}
\mathbf{II}
&\le - \theta h^2 (\psi_k - \psi_{k-1})
     + \gamma \left[ \gamma  - \frac{1}{3}\alpha h\left( \frac{1 + \frac{2}{3}\alpha h}{1+\alpha h} \right) \right] h^4 \| \nabla f(x_k) \|^2\\
&\quad - \frac{\mu h^2}{2}\left[ \gamma  - \frac{1}{3}\alpha h\left( \frac{1 + \frac{2}{3}\alpha h}{1+\alpha h} \right) \right] \| y_{k+1} - y_k \|^2\\
&\le - \theta h^2 (\psi_k - \psi_{k-1})
     + \gamma \left[ \gamma  - \frac{1}{3}\alpha h\left( \frac{1 + \frac{2}{3}\alpha h}{1+\alpha h} \right) \right] h^4 \| \nabla f(x_k) \|^2 \\
&\quad - \frac{\mu\gamma h^2}{2(1+\alpha h)}\| y_{k+1} - y_k \|^2,    
\end{align*}
in view of $\gamma\in[1,2]$ and $\alpha h \le 3$. Whence, we deduce that
\begin{align*}
\left( 1 + \frac{\frac{2}{3}\alpha h}{1 + \frac{2}{3}\alpha h} \right) E_{k+1} - E_k
&\le \frac{\frac{2}{27}\alpha h^3 (1+2\alpha h) \left(\alpha^2 - \frac{9}{2}\mu\gamma \right)}{\left( 1+\frac{2}{3}\alpha h \right)(1+\alpha h)}\| x_k - x^* \|^2 \\
&\quad - \frac{\alpha^2 h^2}{18}\| x_{k+1}-x_k \|^2 
       - R h^4 \| \nabla f(x_k) \|^2 \\ 
&\quad  + \frac{\frac{1}{9} h^2 \left( \alpha^2 - \frac{9}{2}\mu\gamma \right)}{1 + \alpha h } \| y_{k+1} - y_k \|^2,  
\end{align*}
where
\begin{align*}
R
&= \frac{\gamma^2}{2} + \gamma - \frac{ \frac{1}{3}\alpha h }{ 1 + \frac{2}{3}\alpha h } 
  - \left( \gamma - \frac{\frac{1}{3}\alpha h}{1 + \frac{2}{3}\alpha h} \right)\frac{\frac{2}{9}\alpha^2 h^2}{1 + \alpha h } \\
&\quad + \left[ \gamma  - \frac{1}{3}\alpha h\left( \frac{1 + \frac{2}{3}\alpha h}{1+\alpha h} \right) \right] \frac{ \frac{1}{3}\alpha h }{1+\frac{2}{3}\alpha h} 
 - \gamma \left[ \gamma  - \frac{1}{3}\alpha h\left( \frac{1 + \frac{2}{3}\alpha h}{1+\alpha h} \right) \right] \\
&= \frac{\gamma}{2}(2-\gamma) + \frac{ \frac{1}{3}\alpha h (\gamma-1) }{ 1 + \frac{2}{3}\alpha h }
  + \frac{\gamma}{3}\alpha h\left( \frac{1 + \frac{2}{3}\alpha h}{1+\alpha h} \right) \\
&\quad - \left( \gamma - \frac{\frac{1}{3}\alpha h}{1 + \frac{2}{3}\alpha h} \right)\frac{\frac{2}{9}\alpha^2 h^2}{1 + \alpha h }
  - \frac{ \frac{1}{9}\alpha^2 h^2 }{1+\alpha h} \\
&= \frac{\gamma}{2}(2-\gamma) + \frac{ \frac{1}{3}\alpha h (\gamma-1) }{ 1 + \frac{2}{3}\alpha h }
  + \underbrace{ \frac{ \frac{1}{3}\alpha h \left( \gamma - \frac{1}{3}\alpha h \right) }{1+\alpha h} }_{\ge 0} 
  + \underbrace{ \frac{ \frac{2}{27}\alpha^3 h^3 }{\left( 1+\frac{2}{3}\alpha h \right)(1+\alpha h)} }_{\ge 0} \\
&\ge \frac{\gamma}{2}(2-\gamma) + \frac{ \frac{1}{3}\alpha h (\gamma-1) }{ 1 + \frac{2}{3}\alpha h },   
\end{align*}
in view of $\gamma\in[1,2]$ and $\alpha h \le 3$. Then, we obtain 
$$ R\ge \frac{\gamma}{2}(2-\gamma) + \frac{ \frac{1}{3}\alpha h (\gamma-1) }{ 1 + \frac{2}{3}\alpha h } = \frac{ \frac{\gamma}{2}(2-\gamma) + \frac{1}{3}\alpha h \left[ -\left( \gamma - \frac{3}{2} \right)^2 + \frac{5}{4} \right] }{ 1 + \frac{2}{3}\alpha h }
\stackrel{\gamma\in[1,2]}{\ge} \frac{ \frac{1}{3}\alpha h }{ 1 + \frac{2}{3}\alpha h }.  $$
and the desired result follows.
\end{proof}

\section{Proof of Lemma \ref{Lem: E_k_bound_QG}}
\label{Proof: Lem: E_k_bound_QG}

\begin{proof}
By definition of $\phi_k$, we have
\begin{align*}
\phi_{k+1} 
&= (x_{k+1}-y_{k+1}) + \xi h (x_{k+1}-x^*) \\
&= (x_{k+1}-x_k) + \xi h(x_{k+1}-x^*) + h^2 \nabla f(x_k)\\
&= (1+\xi h)(x_{k+1}-x_k)
  + \xi h(x_k-x^*) + h^2 \nabla f(x_k), 
\end{align*}
and using \eqref{E: iterates},
\begin{align*}
\phi_{k+1}-\phi_k
&= -(\alpha-\xi) h (x_{k+1}-x_k) - \gamma h^2 \nabla f(x_k) \\
&= - \frac{(\alpha-\xi) h}{1+2\xi h}\phi_{k+1}
   - \frac{\xi(\alpha-\xi)h^2}{1+2\xi h}(x_{k+1}-x_k)
   + \frac{\xi(\alpha-\xi)h^2}{1+2\xi h}(x_k-x^*) \\
&\quad -\left(\gamma - \frac{(\alpha-\xi) h}{1+2\xi h} \right) h^2 \nabla f(x_k).
\end{align*}
As a result,
\begin{align*}
\| \phi_{k+1}-\phi_k \|^2
&= (\alpha-\xi)^2 h^2\| x_{k+1}-x_k \|^2 
  + \gamma^2 h^4 \| \nabla f(x_k) \|^2 \\
&\quad + 2(\alpha-\xi)\gamma h^3 \langle \nabla f(x_k), x_{k+1} - x_k \rangle,
\end{align*}
and
\begin{align*}
\langle \phi_{k+1}-\phi_k, \phi_{k+1} \rangle
&= - \frac{(\alpha-\xi) h}{1+2\xi h} \| \phi_{k+1} \|^2
   - \frac{\xi(\alpha-\xi) h^2 (1+\xi h )}{1+2\xi h}\| x_{k+1}-x_k \|^2 \\
&\quad + \frac{\xi(\alpha-\xi) h^2}{1+2\xi h}\langle x_{k+1}-x_k, x_k - x^* \rangle \\ 
&\quad - \left[ \gamma (1+\xi h) - \frac{(\alpha-\xi) h}{1+2\xi h} \right] h^2 \langle \nabla f(x_k),x_{k+1}-x_k \rangle\\
&\quad + \frac{\xi^2(\alpha-\xi) h^3}{1+2\xi h}\| x_k - x^* \|^2 
   - \left(\gamma - \frac{2(\alpha-\xi) h}{1+2\xi h} \right) \xi h^3\langle \nabla f(x_k), x_k-x^* \rangle \\
&\quad -\left(\gamma - \frac{(\alpha-\xi) h}{1+2\xi h} \right) h^4 \| \nabla f(x_k) \|^2.
\end{align*}
Since $\frac{1}{2}\| \phi_{k+1} \|^2 - \frac{1}{2}\| \phi_k \|^2 = \langle \phi_{k+1}-\phi_k, \phi_{k+1} \rangle - \frac{1}{2}\| \phi_{k+1} - \phi_k \|^2$, we have
\begin{equation}\label{E: phi_diff_QG}
\begin{aligned}
&\quad \left( 1 + \frac{2(\alpha-\xi) h}{1+2\xi h} \right)\frac{1}{2}\| \phi_{k+1} \|^2 - \frac{1}{2}\| \phi_k \|^2 \\
&= - \frac{\xi(\alpha-\xi) h^2 (1+\xi h )}{1+2\xi h}\| x_{k+1}-x_k \|^2 
   - \frac{(\alpha-\xi)^2 h^2}{2} \| x_{k+1}-x_k \|^2 \\ 
&\quad + \underbrace{ \frac{\xi(\alpha-\xi) h^2}{1+2\xi h} \langle x_{k+1}-x_k, x_k - x^* \rangle }  
  - \left[ \gamma (1+\alpha h) - \frac{(\alpha-\xi) h}{1+2\xi h} \right] h^2 \langle \nabla f(x_k),x_{k+1}-x_k \rangle\\
&\quad + \frac{\xi^2(\alpha-\xi) h^3}{1+2\xi h}\| x_k - x^* \|^2 
   - \left(\gamma - \frac{2(\alpha-\xi) h}{1+2\xi h} \right) \xi h^3\langle \nabla f(x_k), x_k-x^* \rangle \\
&\quad -\left( \frac{\gamma^2}{2} + \gamma - \frac{(\alpha-\xi) h}{1+2\xi h} \right) h^4 \| \nabla f(x_k) \|^2.
\end{aligned}
\end{equation}
Likewise, using \eqref{E: iterates} and recalling that
$$ y_{k+1}-x^* = x_k - x^* - h^2\nabla f(x_k), $$
we have
\begin{align*}
y_{k+1} - y_k
&= (1+\alpha h)(x_{k+1}-x_k) + \gamma h^2 \nabla f(x_k) \\
&= -\frac{(\alpha-\xi) h}{1+2\xi h}(y_{k+1}-x^*)
   + (1+\alpha h)(x_{k+1}-x_k) 
   + \frac{(\alpha-\xi) h}{1+2\xi h}(x_k - x^*) \\
&\quad + \left( \gamma - \frac{(\alpha-\xi) h}{1+2\xi h} \right) h^2 \nabla f(x_k), 
\end{align*}
so that
\begin{align*}
\| y_{k+1} - y_k \|^2
&= (1+\alpha h)^2 \| x_{k+1}-x_k \|^2 + \gamma^2 h^4 \| \nabla f(x_k) \|^2 \\
&\quad + 2(1+\alpha h)\gamma h^2\langle \nabla f(x_k), x_{k+1}-x_k \rangle,
\end{align*}
and
\begin{align*}
&\quad \langle y_{k+1}-y_k, y_{k+1}-x^* \rangle \\  
&= -\frac{(\alpha-\xi) h}{1+2\xi h}\| y_{k+1}-x^* \|^2
  + (1+\alpha h)\langle x_{k+1}-x_k, x_k - x^* \rangle\\
&\quad - (1+\alpha h)h^2 \langle \nabla f(x_k), x_{k+1}-x_k \rangle
  + \frac{(\alpha-\xi) h}{1+2\xi h} \| x_k - x^* \|^2 \\
&\quad + \left( \gamma - \frac{2(\alpha-\xi) h}{1+2\xi h} \right) h^2 \langle \nabla f(x_k), x_k - x^* \rangle   
  - \left( \gamma - \frac{(\alpha-\xi) h}{1+2\xi h} \right) h^4 \| \nabla f(x_k) \|^2.
\end{align*}
Since $\frac{1}{2}\| y_{k+1}-x^* \|^2 - \frac{1}{2}\| y_k-x^* \|^2 = \langle y_{k+1}-y_k, y_{k+1}-x^* \rangle - \frac{1}{2}\| y_{k+1} - y_k \|^2$, we have
\begin{align*}
&\quad \left( 1 + \frac{2(\alpha-\xi) h}{1+2\xi h} \right)\frac{1}{2}\| y_{k+1}-x^* \|^2 - \frac{1}{2}\| y_k-x^* \|^2 \\
&= \underbrace{ (1+\alpha h)\langle x_{k+1}-x_k, x_k - x^* \rangle }
   - \frac{1}{2}(1+\alpha h)^2 \| x_{k+1}-x_k \|^2 \\
&\quad - (1+\alpha h)(1+\gamma)h^2 \langle \nabla f(x_k), x_{k+1}-x_k \rangle
  + \frac{(\alpha-\xi) h}{1+2\xi h} \| x_k - x^* \|^2 \\
&\quad + \left( \gamma - \frac{2(\alpha-\xi) h}{1+2\xi h} \right) h^2 \langle \nabla f(x_k), x_k - x^* \rangle 
  - \left( \frac{\gamma^2}{2} + \gamma - \frac{(\alpha-\xi) h}{1+2\xi h} \right) h^4 \| \nabla f(x_k) \|^2. 
\end{align*}
To cancel the term $\langle x_{k+1}-x_k, x_k - x^* \rangle$ in \eqref{E: phi_diff_QG}, we multiply this inequality by
$$ \eta = \frac{\xi(\alpha-\xi) h^2}{(1+\alpha h)(1+2\xi h)} $$
and subtract it from \eqref{E: phi_diff_QG}, to get
\begin{align*}
&\quad \left( 1 + \frac{2(\alpha-\xi) h}{1+2\xi h} \right) \left( \frac{\| \phi_{k+1} \|^2}{2} - \frac{\eta}{2}\| y_{k+1}-x^* \|^2 \right) 
  - \left( \frac{\| \phi_k \|^2}{2} - \frac{\eta}{2}\| y_k-x^* \|^2 \right)\\
&= \underbrace{ - \frac{\xi(\alpha-\xi) h^2 [1+(2\xi-\alpha) h ]}{2(1+2\xi h)} \| x_{k+1}-x_k \|^2 
   - \frac{(\alpha-\xi)^2 h^2}{2}\| x_{k+1}-x_k \|^2 }_{\le 0}\\
&\quad \underbrace{ - \left[ \gamma - \frac{(\alpha-\xi)h \left( 1 + \xi h (1+\gamma) \right) }{(1+\alpha h)( 1+2\xi h )} \right] (1+\alpha h) h^2 \langle \nabla f(x_k),x_{k+1}-x_k \rangle }_{\mathbf{I}}\\
&\quad + \underbrace{ \frac{\xi(\alpha-\xi)^2 h^3 \left[ \frac{\xi}{\alpha-\xi}(1+\alpha h)(1+2\xi h) - 1 \right] }{(1+\alpha h)\left(1+2\xi h\right)^2}\| x_k - x^* \|^2 }_{\mathbf{II}} \\ 
&\quad \underbrace{ - \left(\gamma - \frac{2(\alpha-\xi) h}{1+2\xi h} \right) \left[ \frac{ (1+2\alpha h)( 1 + \xi h ) }{(1+\alpha h)( 1 + 2\xi h )} \right] \xi h^3\langle \nabla f(x_k), x_k-x^* \rangle }_{\mathbf{III}} \\
&\quad -\left( \frac{\gamma^2}{2} + \gamma - \frac{(\alpha-\xi) h}{1+2\xi h} \right)(1-\eta) h^4 \| \nabla f(x_k) \|^2.
\end{align*}
Recalling that
$$ \theta = \gamma - \frac{(\alpha-\xi)h \left( 1 + \xi h (1+\gamma) \right) }{(1+\alpha h)( 1+2\xi h )}, $$
we use Lemma \ref{Lem: func_bound}, to obtain
\begin{align*}
\mathbf{I}
&\le -\theta h^2 (\psi_k-\psi_{k-1}) + \gamma \left[ \gamma - \frac{(\alpha-\xi)h \left( 1 + \xi h (1+\gamma) \right) }{(1+\alpha h)( 1+2\xi h )} \right] h^4 \| \nabla f(x_k) \|^2 \\
&= -\theta h^2 (\psi_k-\psi_{k-1}) 
   + \gamma \left[ \gamma - \frac{\xi(\alpha-\xi)h^2(1+\gamma)}{(1+\alpha h)( 1+2\xi h )} \right] h^4 \| \nabla f(x_k) \|^2 \\
&\quad - \frac{ \gamma(\alpha-\xi) h }{(1+\alpha h)( 1+2\xi h )}h^4 \| \nabla f(x_k) \|^2 \\
&= -\theta h^2 (\psi_k-\psi_{k-1}) 
   + \gamma\left[ (1-\eta)\gamma - \eta \right]h^4 \| \nabla f(x_k) \|^2 
   - \frac{ \gamma(\alpha-\xi) h^5 \| \nabla f(x_k) \|^2 }{(1+\alpha h) ( 1+2\xi h )}.
\end{align*} 
Likewise, we deduce that
\begin{align*}
\mathbf{II} 
&= \frac{\xi(\alpha-\xi)^2 h^3}{(1+\alpha h)\left(1+2\xi h\right)^2} \left[ \frac{2\xi-\alpha}{\alpha-\xi} + \frac{\xi}{\alpha-\xi}(\alpha+2\xi) h + \frac{\xi}{\alpha-\xi} 2\alpha\xi h^2 \right] \| x_k - x^* \|^2 \\
&= \frac{\xi(\alpha-\xi)(2\xi-\alpha) h^3}{(1+\alpha h)\left(1+2\xi h\right)^2}\left[ 1 + \frac{\xi}{2\xi-\alpha}(\alpha+2\xi) h + \frac{\xi}{2\xi-\alpha} 2\alpha\xi h^2  \right]\| x_k - x^* \|^2 \\ 
&\stackrel{\eqref{E: ineqn_frac_quad}}{\le} 
  \frac{\xi(\alpha-\xi)(2\xi-\alpha) h^3}{(1+\alpha h)\left(1+2\xi h\right)}\left( 1 + \frac{\xi(3\alpha-2\xi)h}{2\xi-\alpha} \right) \| x_k - x^* \|^2,
\end{align*}
and in view of convexity,
\begin{align*}
\mathbf{III}
&\le - \left(\gamma - \frac{2(\alpha-\xi) h}{1+2\xi h} \right) \left[ \frac{ (1+2\alpha h)( 1 + \xi h ) }{(1+\alpha h)( 1 + 2\xi h )} \right] \xi h^3 \left( f(x_k) - f^* \right) \\ 
&= - \left[ \frac{\xi}{2(\alpha-\xi)}(1+2\alpha h)(1+\xi h)\gamma - \xi h \frac{(1+2\alpha h)( 1 + \xi h )}{ 1 + 2\xi h } \right] \frac{ 2(\alpha-\xi) h^3 \left( f(x_k) - f^* \right) }{(1+\alpha h)(1+2\xi h)} \\
&\stackrel{\eqref{E: ineqn_frac}}{\le} 
  - \left[ \frac{\xi}{2(\alpha-\xi)} (1+2\alpha h)(1+\xi h)\gamma - \xi h \left(1+(2\alpha-\xi) h\right) \right] \frac{ 2(\alpha-\xi) h^3 \left( f(x_k) - f^* \right) }{(1+\alpha h)(1+2\xi h)}\\
&\le
  - \underbrace{ \left[ \frac{\xi}{2(\alpha-\xi)}\left(1+(2\alpha+\xi) h\right)\gamma - \xi h \right] }_{\mathbf{IV}} \frac{ 2(\alpha-\xi) h^3 \left( f(x_k) - f^* \right) }{(1+\alpha h)(1+2\xi h)}, 
\end{align*}
where we used $\gamma\ge 1$ and $\frac{\xi}{2(\alpha-\xi)}\ge 1$ in the last inequality. Similarly, using $\gamma\ge 1$ and $\frac{\xi}{2(\alpha-\xi)}\ge 1$, we have 
\begin{align*}
\mathbf{IV}
&= (1+\alpha h)\gamma + \frac{3\xi-2\alpha}{2(\alpha-\xi)} \gamma 
  + \frac{\xi}{2(\alpha-\xi)}\gamma \alpha h  
  + \underbrace{ \frac{\xi}{2(\alpha-\xi)}\gamma(\alpha+\xi)h - (\gamma \alpha+\xi)h }_{\ge 0} \\
&\ge (1+\alpha h)\gamma + \frac{3\xi-2\alpha}{2(\alpha-\xi)} \gamma \left( 1 + \frac{\alpha\xi h}{3\xi-2\alpha} \right) \\
&\ge (1+\alpha h)\gamma + \frac{3\xi-2\alpha}{2(\alpha-\xi)} \gamma \left( 1 + \frac{\xi(3\alpha-2\xi)h}{2\xi-\alpha} \right) \quad \left( \text{since } \xi = \frac{2+\sqrt{2}}{3+\sqrt{2}}\alpha \right) \\
&= (1+\alpha h)\gamma + \frac{\gamma}{\sqrt{2}} \left( 1 + \frac{\xi(3\alpha-2\xi)h}{2\xi-\alpha} \right). 
\end{align*}
Therefore,
\begin{align*}
\mathbf{III}
&\le - \left[ (1+\alpha h)\gamma + \frac{\gamma}{\sqrt{2}}\left( 1 + \frac{\xi(3\alpha-2\xi)h}{2\xi-\alpha} \right) \right] \frac{ 2(\alpha-\xi) h^3 \left( f(x_k) - f^* \right) }{(1+\alpha h)(1+2\xi h)} \\
&= - \frac{ 2\gamma(\alpha-\xi) h^3 }{1+2\xi h}\left(\psi_k + \frac{h^2}{2}\| \nabla f(x_k) \|^2\right)
    - \frac{ \sqrt{2}\gamma(\alpha-\xi) h^3 \left( 1 + \frac{\xi(3\alpha-2\xi)h}{2\xi-\alpha} \right) }{(1+\alpha h)\left(1+2\xi h\right)} \left( f(x_k) - f^* \right) \\
&\stackrel{\theta\le\gamma}{\le} 
    - \left( \frac{ 2(\alpha-\xi) h }{1+2\xi h} \right)\theta h^2 \psi_k
    - \frac{ \gamma(\alpha-\xi) h^5 }{1+2\xi h} \| \nabla f(x_k) \|^2 \\
&\quad - \frac{ \frac{\mu\gamma}{\sqrt{2}}(\alpha-\xi) h^3 \left( 1 + \frac{\xi(3\alpha-2\xi)h}{2\xi-\alpha} \right) }{(1+\alpha h)\left(1+2\xi h\right)} \| x_k - x^* \|^2, 
\end{align*}
where we applied the quadratic growth condition in the last inequality. Whence,
\begin{align*}
\left( 1 + \frac{ 2(\alpha-\xi) h }{1+2\xi h} \right)E_{k+1} - E_k
&\le \frac{ (\alpha-\xi) h^3 \left( 1 + \frac{\xi(3\alpha-2\xi)h}{2\xi-\alpha} \right) \left[ \xi(2\xi-\alpha) - \frac{\mu\gamma}{\sqrt{2}} \right]}{(1+\alpha h)\left(1+2\xi h\right)}\| x_k - x^* \|^2 \\ 
&\quad - R h^4 \| \nabla f(x_k) \|^2,
\end{align*}
where
\begin{align*}
R 
&= \frac{\gamma}{2}(2-\gamma) + \frac{\eta\gamma^2}{2} + \frac{ (\alpha-\xi) h (\gamma-1+\eta) }{ 1 + 2\xi h } + \frac{ \gamma(\alpha-\xi) h }{(1+\alpha h)( 1 + 2\xi h )} \\
&\ge \frac{\gamma}{2}(2-\gamma) + \frac{ (\alpha-\xi) h (\gamma-1) }{ 1 + 2\xi h } \\
&= \frac{ \frac{\gamma}{2}(2-\gamma) + (\alpha-\xi) h \left[ (2+\sqrt{2})\gamma(2-\gamma) + (\gamma-1) \right]  }{ 1 + 2\xi h } \quad \left(\text{since }\xi = (2+ \sqrt{2}) (\alpha-\xi) \right)\\
&\stackrel{\gamma\in[1,2]}\ge 
  \frac{ (\alpha-\xi) h}{ 1 + 2\xi h }.   
\end{align*}
The result follows after reorganizing the terms.
\end{proof}

\end{appendices}

\end{document}